\documentclass[12pt,a4paper]{amsart}
\usepackage[utf8]{inputenc}	
\usepackage[T1]{fontenc}
\usepackage[in,headings]{fullpage}
\usepackage{amsbsy, amscd, amsfonts, amsmath, amsrefs, amssymb, amsthm}
\usepackage{mathrsfs}
\usepackage{graphicx}
\usepackage{float}
\usepackage{color}   
\usepackage[colorlinks=true,linkcolor=blue,citecolor=blue,urlcolor=blue]{hyperref}
\usepackage{setspace}

\newtheorem{theorem}{Theorem}

\newtheorem{proposition}[theorem]{Proposition}

\theoremstyle{definition}
\newtheorem{definition}[theorem]{Definition}
\newtheorem{remark}[theorem]{Remark}

\theoremstyle{remark}

\newcommand{\Z}{\mathbf{Z}}

\newcommand{\R}{\mathbf{R}}
\newcommand{\N}{\mathbf{N}}

\newcommand{\cpx}[1]{\|#1\|}
\newcommand{\dft}{\delta}
\newcommand{\st}{\mathrm{st}}

\begin{document}

\title[Arithmetic Compact Sets]
{Complexity of Natural numbers and arithmetic Compact sets.}
\author[Arias de Reyna]{J. Arias de Reyna}
\address{%
Universidad de Sevilla \\ 
Facultad de Matem\'aticas \\ 
c/Tarfia, sn \\ 
41012-Sevilla \\ 
Spain.} 




\email{arias@us.es, ariasdereyna1947@gmail.com}

\date{\today, \texttt{201-TheCompact-v3.tex}}

\begin{abstract}
The complexity $\Vert n\Vert$ of a natural number is the least number of $1$ needed to represent $n$ using the 5 symbols $(, ), *, +, 1$. A natural number $n$ is called stable is $\cpx{3^kn}=\cpx{n}+3k$. For each natural number $n$, the number $3^an$ is stable for some $a\ge0$, and we define the stable complexity of $n$ as $\cpx{n}_\st=\cpx{3^an}-3a$. We show that the closure of the set of all fractions $n/3^{\lfloor \cpx{n}_\st/3\rfloor}$ has remarkable properties; self-similarity $3K'''=K$, well-ordered, and certain arithmetical properties.  We pose the question about the unicity of this compact. This question raises some problems about the complexity of natural numbers that we are unable to answer.
\end{abstract}

\maketitle

\section{Introduction}

The complexity $\cpx{n}$ of a natural number $n$ is the least number of $1$'s needed to write an expression for $n$ in the language using only the addition $+$, the product $*$, the unit $1$  and parentheses as only symbols. So $\cpx{9}=6$ since $9=(1+1+1)*(1+1+1)$ and there is no other expression for $9$ with less than $6$ unit symbols. Richard Guy popularised several problems on this concept, such as asking if $\cpx{2^n}=2n$ for all $n\ge1$.  

Frequently we have $\cpx{3n}=\cpx{n}+3$. Therefore, it is natural to think about the fractions $n/3^{\lfloor\cpx{n}/3\rfloor}$ that remain invariant when we change $n$ by $3n$ and $\cpx{3n}=\cpx{n}+3$. 
A number is called stable when $\cpx{n3^k}=\cpx{n}+3k$. And for each $m$ there is some $a$ with $n=m3^a$ stable. In \cite{Arias} I made some conjectures about these fractions. These conjectures have been solved by H. Altman \cite{paper1}, \cite{paperwo} and H. Altman, J. Arias de Reyna \cite{AA}. In this paper, I give a new interpretation of these conjectures. All these conjectures are summarised in the existence of a surprising object, which I call an \emph{arithmetic compact set}.

As a young man, I lived isolated from other mathematicians and was always eager to learn about open problems. When in 2000, as a professional mathematician, I published a paper in the `Gaceta' in Spanish full of conjectures, I was thinking of a young looking for a topic to research. I wanted to give a fruitful topic. I was not interested in asking accurate questions, but leave to this young the formulation of the correct statements. I knew that what I wrote was not exactly true, but I was sure that I was directing this young audience to a rich topic. 

These objectives were amply met when Harry Altman and Josh Zelinsky submitted their paper to \emph{Integers} in late 2011. After acting as a referee for \cite{paper1}, I ended up collaborating with H.~Altman on \cite{AA} in which we reformulate and resolve the rest of my conjectures,   and also in some other subsequent work \cite{hyperplanes}, \cite{deg1}, on the complexity.  In the first Section of this paper I resume the state of the conjectures. Most of them are true after some reformulation. 

The main purpose of this paper is to show how these conjectures are equivalent to the existence of a compact set $K$ of rational numbers with some remarkable properties. Properties of self-similarity and also of an arithmetical nature. The self-similarity is synthesized in the equation $3K'''=K$. The compact $K$ is also well-ordered by the reverse usual order. For each natural number $m$ with $3\nmid m$ there is a non-negative integer $\kappa(m)$ such that 
\[K=\{0\}\cup\Bigl\{\frac{m}{3^k}\colon 3\nmid m,  k\ge \kappa(m)\Bigr\}.\]
This are not all properties of this compact set. 

We have another main example of an arithmetical compact, verifying $2K'=K$. Therefore, in Section 3 we study a general concept. Two examples may not be enough to justify this, but will certainly avoid repetition when I explain the second example in a later paper. They also make the role of each conjecture in \cite{Arias} clearer. 

Section \ref{existence} is devoted to proving the existence of the arithmetical compact $K$. 
This consists in translating the results of \cite{paper1}, \cite{paperwo} and \cite{AA} where my conjectures are refined and proved, into the language of \cite{Arias}. 

When I saw the existence of the compact, I saw that its elements appeared in an almost inevitable way, suggesting its uniqueness. In Section \ref{unicity} I ask about this: is $K$ unique? Any derivative set $K^{(n)}$ is also an arithmetical compact. My question is after considering all of these as a unique solution. The answer is probably no. This is due to the construction by H. Altman of another arithmetical compact $H$. But the problem to prove $H\ne K$ is not trivial. In the Appendix we give some of the first elements in $K$. All these elements are also the first elements in $H$. In this Section we consider possible candidates of fractions $q\in K\smallsetminus H$.  But the proof is beyond our capacity to compute at this time.

\section{State of the conjectures}
The main contribution of \cite{paper1} is the definition of the \emph{defect}
$\dft(n):=\cpx{n}-3\log_3n$, 
and proving that for any $x>0$ the set 
\[A_x=\{n\in\N\colon \dft(n)<x\}\]
has a simple alternate description that can be computed efficiently. 

There is a change of language between the papers
\cite{paper1}, \cite{paperwo}, \cite{AA} and  \cite{Arias}. To say that 
$\mathscr{D}^u_\st[\alpha]=\dft_\st(n)$ means that $\dft_\st(n)=\cpx{n}_\st-3\log_3n$ and $\cpx{n}_\st=3\ell+u$ for some $u\in\{0,1,2\}$. Therefore $\frac{n}{3^\ell}=S_u[\alpha]$. And then
\[u-3\log_3S_u[\alpha]=u-3\log_3\frac{n}{3^\ell}=3\ell+u-3\log_3n=
\dft_\st(n)=\mathscr{D}^u_\st[\alpha].\]
This relationship will make translation easy in the following. 

\subsection{Conjecture 1} (True).  A number $n$ is called \emph{stable} if $\cpx{3^kn}=\cpx{n}+3k$ for all $k\ge0$. The conjecture says that for any natural number $n$ there is $K$ such that $n3^K$ is stable. This is proved in \cite{paper1}*{Thm.~13}. The \emph{stable complexity} of $n$ is then defined by $\cpx{n}_\st:=\cpx{3^K n}-3K$. 

\subsection{Conjecture 2} (True after reformulation). The Conjecture states that 
$\cpx{p(q3^j+1)}=\cpx{p}+\cpx{q}+3j+1$. This is true if $pq$ is stable and $\cpx{pq}=\cpx{p}+\cpx{q}$. This is proved in \cite{AA}*{Thm.~1.4}. This condition was not mentioned in the original statement of Conjecture 2. They are needed, in \cite{AA} the counterexample $p=2$, $q=1094$ is mentioned, where 
\[\cpx{2\cdot1094}=\cpx{3^7+1}=22<\cpx{2}+\cpx{1094}=24.\]
And we have 
$\cpx{2(1094\cdot 3^k+1)}=\cpx{2188\cdot3^k+2}\le 24+3k,$
instead of $\cpx{2(1094\cdot 3^k+1)}=25+3k$.

In \cite{AA}*{Thm.~1.18} another version of the conjecture is proved with $\cpx{pq}+1=\cpx{p}+\cpx{q}$. 

\subsection{Conjecture 3} (As stated, it is false). Three sets are defined in Conjectures 5, 6 and 7. They are well ordered, which is proved in \cite{paperwo}*{Thm.~A.5}, but the main assertion of Conjecture 3, that the (gratest) numbers of complexity $3n+u$  are the first natural numbers contained in $3^n\times\text{these sequences}$, is not true, because these numbers are not all stable. We may reformulate this, as the stable numbers of complexity $3n+u$ are equal to the stable numbers in $3^n\times\text{these sequences}$. But then this is almost tautological.

With good will, we could say that the part of the good order in Conjecture 3 is true but not the rest.

\subsection{Conjectures 4, 5, 6, 7} (True).  Conjecture 4 is true by definition if we take the definition in Conjectures 5, 6 and 7 of the sequences. The three defined sets are
\begin{equation}\label{def:S}
S_u:= \Bigl\{\frac{n}{3^k}\colon n \text{ stable and } \cpx{n}=3k+u\Bigr\},\qquad u\in\{0,1,2\}.
\end{equation}
The conjectures (or is this in Conjecture 3?) state that they are well ordered by the reverse usual order of $\R$. 

As said before, 5, 6 and 7 are proved in \cite{paperwo}*{Thm.~A.5}.

Summary: Conjecture 3 says that there are some sequences with two properties: good order and representation of certain numbers, the second part is not true. Conjectures 5, 6 and 7 define some sequences and say that they coincide with those of Conjecture 3. The property of good order predicted in Conjecture 3 is fulfilled in the sequences defined in Conjectures 5, 6 and 7.

\subsection{Conjecture 8} (True). The sets $S_u$ are well ordered by $x\preccurlyeq y$ if and only if $x>y$. For a well-ordered set $L$ and an ordinal $\alpha$ we denote by $L[\alpha]$ the element of $L$ in position $\alpha$. With these notations, Conjecture 8, asserts that
\begin{equation}\label{conj8}
\lim_{n\to\infty}S_0[\omega\alpha+n]=S_2[\alpha]/3,\quad \lim_{n\to\infty}S_1[\omega\alpha+n]=S_0[\alpha],\quad \lim_{n\to\infty}S_2[\omega\alpha+n]=S_1[\alpha].
\end{equation}
This is proved in \cite{AA}*{Cor.~4.17}, (the change of language offers no serious difficulty).

\subsection{Conjectures 9, 10 and 11} (True after reformulation). Theorem 4.19 in \cite{AA} gives the true version of these conjectures. We may state it in our language as 
\begin{theorem}\label{limitpoints}
For any ordinal $\beta<\omega^\omega$ let $S_u[\beta]=\frac{n}{3^k}$ with $\cpx{n}_\st=3k+u$, the two sets 
\[\{S_{u+1}[\omega\beta+r]\colon r\in\Z_{\ge0}\},\quad\Bigl\{\frac{b(a 3^r+1)}{3^{r+k+\varepsilon}}\colon n=ab,\quad \cpx{n}_{\st}=\cpx{a}_\st+\cpx{b}_\st, \quad r\in\Z_{\ge0}\Bigr\}\]
have a finite symmetric difference.
(Understanding $u+1\bmod 3$ and $\varepsilon=0$ for $u=0$, $1$ and $\varepsilon=1$ for $u=2$). 
\end{theorem}
The conditions on $n=ab$, $a$ and $b$ in the second set were not correct in the statement of Conjectures 9, 10 and 11. Conjectures 9, 10 and 11 speak of the difference between the two sets consists of the fractions corresponding to powers of 2. This is not true, there are powers of two, but for example $S_0(\omega2) =\frac{1280}{2187}$ is a counterexample appearing in the Gaceta paper!

I will return to this conjecture and these sporadic numbers later. 
I add a table summarizing the situation

\begin{center}
\begin{tabular}{|c|l|c|} \hline
\multicolumn{3}{|c|}{\bfseries Gaceta conjectures}\\ \hline
Conjecture & Truth value & Solution found in \\ 
\hline \noalign{\vspace{0.5cm}}\hline
1 & True & \cite{paper1}*{Thm.~13}\\ \hline
2 & True (re-stated) & \cite{AA}*{Thm.~1.4 and 1.18}\\ \hline
3(wo), 4, 5, 6, 7 & True & \cite{paperwo}\\ \hline
3(r) & False &  \\ \hline
8 & True & \cite{AA}*{Cor.~4.17}\\ \hline
9, 10, 11  main part& True (re-stated) & \cite{AA}*{Thm.~4.19} \\ \hline
9, 10, 11 sporadics & False & \\
\hline
\end{tabular}
\end{center}

\section{Arithmetic self-similar compact sets}

I have only two examples for the concept that I present here. My objective is to present one of these examples and ask a question about the unicity of this example; the formal definition seems to be the easiest way to formulate this question.

For any set $A\subset \R$ we denote by $A'$ the derived set of $A$, the set of accumulation points of $A$. Also, for any nonnegative integer $n$, let $ A^{(n)}$ the $n$-th derivative set, with $A^{(0)}:= A$.

We will consider subsets of the set of real numbers $\R$, well ordered by $\preccurlyeq$,  the \emph{reverse usual order}. That is,  $x\preccurlyeq y$ if and only if $x\ge y$.

A well-ordered set is isomorphic to an ordinal. If $K$ is well ordered with the order type $\omega^\omega+1$, and $0\le\alpha\le\omega^\omega$, we denote by   $K[\alpha]$ the corresponding element in $K$. Hence, in our examples we will have $K\subset[K[\omega^\omega], K[0]]$. 

\begin{definition}
Let $\rho$ be a real number  and $\mu$ a natural number with $\rho>1$.
A compact set $K\subset[0,+\infty)$ will be called a self-similar compact set of ratio $\rho$ and module $\mu$ if it satisfies 
\begin{equation}
\rho K^{(\mu)}=K
\end{equation} and it is well ordered by the reverse usual order $\preccurlyeq$.
\end{definition}
 
\begin{proposition}\label{Kpos}
Let $K$ be a self-similar compact set. Then $K[\omega^\omega]=0$, $K[0]>0$ and $K\subset[0, K[0]]$. 
\end{proposition}
\begin{proof}
Since $K$ is compact, there is some $a>0$ with $K\subset[-a,a]$. Since $K$ is self-similar, we have $K^{(\mu)}=K/\rho$. By successively deriving this relation, by  induction, we get 
$K^{(n\mu)}=K/\rho^n\subset[-a/\rho^n, a/\rho^n]$. Since $\rho>1$, it follows that 
$\bigcap_{n=0}^\infty K^{(n)}\subset\{0\}$. The order type of $K$ is $\omega^\omega+1$, and this implies that $\bigcap_{n=0}^\infty K^{(n)}=K^{(\omega)}=\{K[\omega^\omega]\}$. Therefore, $K[\omega^\omega]=0$. For any other element $y\in K$ we have  $K[0]\preccurlyeq y\preccurlyeq K[\omega^\omega]=0$. Therefore, $0\le y\le K[0]$. We have $K[0]>0$, since $K[0]\ne K[\omega^\omega]$.
\end{proof}
\begin{proposition}\label{P:cuatro}
Let $K$ be a self-similar compact set of modulo $\mu$ and ratio $\rho$, then we have 
\begin{equation}\label{230118-2}
K=\{0\}\cup\bigcup_{n=0}^\infty(T_0\cup T_1\cup\cdots\cup T_{\mu-1})/\rho^n,
\end{equation}
where $T_u=K^{(u)}\smallsetminus K^{(u+1)}$.
\end{proposition}
\begin{proof}
For $n\ge0$ we have $K^{(n)}\supset K^{(n+1)}$, and 
$\bigcap_{n=0}^\infty K^{(n)}=\{0\}$ as we have seen in Proposition \ref{Kpos}. Therefore, if we define $T_u=K^{(u)}\smallsetminus K^{(u+1)}$ for $u\in\Z_{\geq0}$, these sets are disjoints and 
\begin{equation}\label{230118-1}
K\smallsetminus\{0\}=\bigcup_{u=0}^\infty T_u.
\end{equation}
By successive derivation from $K^{(\mu)}=K/\rho$ we get $K^{(u+n\mu)}=K^{(u)}/\rho^n$. Therefore, $T_{u+n\mu}=T_u/\rho^n$. Substituting this into \eqref{230118-1} we get \eqref{230118-2}.
\end{proof}
\begin{proposition}
Let $K$ be a self-similar compact set of module $\mu$ and ratio $\rho$, and $\alpha$ an ordinal with $0\le \alpha\le \omega^{\omega}$ we have
\begin{gather*}
 K^{(n)}[\alpha]=K[\omega^n(1+\alpha)],\qquad n\ge 1.\\
T_0[\alpha]=\begin{cases} K[\alpha] & \text{if $\alpha<\omega$,}\\
K[\alpha+1] & \text{if $\alpha\ge\omega$,}\end{cases}\quad\text{and}\quad T_u[\alpha]=K[\omega^u(\alpha+1)], \quad\text{ for $u\ge1$}.
\end{gather*}
\end{proposition}
\begin{proof}
The set $K'$ is the set of limit points of $K$, which corresponds to the limit ordinals in $\omega^\omega+1$. Ordinals that can be written as $\alpha=\omega(\alpha+1)$. Therefore, 
$K'[\alpha]=K[\omega(1+\alpha)]$. It follows that 
\[K''[\alpha]=K'[\omega(1+\alpha)]=K[\omega(1+\omega(1+\alpha))]=K[\omega^2(1+\alpha)]\]
And by induction, we get $K^{(n)}[\alpha]=K[\omega^n(1+\alpha)]$. 

The elements in $T_0$ correspond to the set of ordinal $0$ and the successors. The application that sends $\alpha$ to $\alpha$ or $\alpha+1$ acoording to $\alpha<\omega$  or $\alpha\ge\omega$ is injective, increasing and takes all successors and 0. This implies the formula for $T_0[\alpha]$ in our Proposition. For $u>0$ the elements in $T_u$ are those contained in the form $K[\omega^u(1+\alpha)]$ that are not limit points of these ones. That is $K[\omega^u(\alpha+1)]$.
\end{proof}
\begin{proposition}
For  integers $0\le\alpha<\omega^\omega$ and $u\in\Z_{\ge0}$ we have 
\[\lim_{n\to\infty}T_u[\omega\alpha+n]=T_{u+1}[\alpha].\]
When $u=\mu-1$ we may write this in the form 
\[\lim_{n\to\infty}T_{\mu-1}[\omega\alpha+n]=\frac{T_{0}[\alpha]}{\rho}.\]
\end{proposition}

\begin{proof}
Since $K$ is a compact, its limits points are also limits points in the topological sense, so that 
\[\lim_{n\to\infty}K[\omega\alpha+n]=K[\omega(\alpha+1)]\quad\text{and}
\lim_{n\to\infty}K[\omega^a(\omega\alpha+n)]=K[\omega^{a+1}(\alpha+1)]\quad \text{for $a\in \N$}\]
Hence, for $u\ge0$ we have 
\[\lim_{n\to\infty} T_u[\omega\alpha+n]=\lim_{n\to\infty}K[\omega^u(\omega\alpha+n+1)]=
K[\omega^{u+1}(\alpha+1)]=T_{u+1}[\alpha].\]
(For $u=0$ and $\alpha=0$ notice that the limit when $n$ or $n+1$ tends to infinity is the same thing). 

Now, if $u=\mu-1$ notice that 
\[T_\mu=K^{(\mu)}\smallsetminus K^{(\mu+1)}=\rho^{-1}(K\smallsetminus K')=\rho^{-1}T_0.\qedhere\]
\end{proof}
We are going to define the \emph{arithmetic} self-similar compact sets of ratio $\rho$ and module $\mu$. It is preferable to give the definition in two steps. First, we define the quasi-arithmetic compact.  In this way we may give the definition of the stability gauge $\cpx{n}_K$, before giving the interesting last definition of arithmetic compact set.

If $K$ is a self-similar compact set and $x\in K$ with $x\ne0$, we have 
\[\cdots  \preccurlyeq x\rho^n \preccurlyeq x\rho^{n-1}\preccurlyeq\cdots\preccurlyeq x\rho \preccurlyeq x\]
Since $K$ is well ordered there is some $n\ge0$ such that $y:=x\rho^n\in K$ but $y\rho\notin K$. 
and $x=y\rho^{-n}$. 
$K$ is quasi-arithmetic when $\rho\ge2$ is a natural number and all these $y=m/\rho^\kappa$ with $\kappa\in\Z_{\ge0}$,  $m\in \N$ and $\rho\nmid m$. More precisely:

\begin{definition}
A self-similar compact set $K$ of module $\mu$ and ratio $\rho$ is called {quasi-arithmetic} if $\rho\ge2$ is a natural number, and there is a function $\kappa\colon\N\smallsetminus (\rho\N)\to \Z_{\ge0}$ such that 
\[K\smallsetminus\{0\}=\Bigl\{\frac{n}{\rho^k}\colon \rho\nmid n \text{ and } k\ge\kappa(n)\Bigr\},\]
\end{definition}
\begin{proposition}\label{P_230127-1}
Let $K$ be a quasi-arithmetic compact set. Assume that $\frac{n}{\rho^p}\in T_v$ and $\frac{n}{\rho^q}\in T_w$. Then $p\mu-v=q\mu-w$.
\end{proposition}
\begin{proof}
Let $m\in\N$ be such that $\rho\nmid m$ and put $k=\kappa(m)$, equation \eqref{230118-2} implies that $\frac{m}{\rho^{k}}\in T_0\cup T_1\cup\cdots\cup T_{\mu-1}$. Since these sets are disjoint, there is a unique $u\in\{0,1,\dots,\mu-1\}$ with $\frac{m}{\rho^{k}}\in T_u$. 
Then $\frac{m}{\rho^{k+\ell}}\in T_{u+\ell\mu}$. Since the $T_v$ are disjoints, these are the unique $T_v$ containing these fractions.

As in the proposition assume that $\frac{n}{\rho^p}\in T_v$ and $\frac{n}{\rho^q}\in T_w$, and let $n=m\rho^r$ with $\rho\nmid m$ and $r\ge0$. With the notation above, we will have 
for some $\ell$ that  $p-r=k+\ell$ and $v=u+\ell\mu$. Therefore, 
\[p\mu-v=(r+k+\ell)\mu-(u+\ell \mu)=(r+k)\mu-u=q\mu-w.\qedhere\]
\end{proof}

\begin{definition}
Let $K$ be a quasi-arithmetic compact set $K$ of module $\mu$ and ratio $\rho$.
For any natural number $n$ we define the \emph{stability gauge} $\cpx{n}_K$  of $n$ by 
\[\cpx{n}_K=\ell\mu+d-u,\text{ when } \frac{n}{\rho^\ell}\in T_u\]
where the displacement $d$ is chosen so that $\cpx{1}_K=0$. 
\end{definition}

Note that by Proposition \ref{P_230127-1},  $\cpx{n}_K$ do not depend on the chosen $\ell$ and $u$ with $\frac{n}{\rho^\ell}\in T_u$. 

\begin{proposition}
For any natural number $n$ we have $\cpx{\rho n}_K=\cpx{n}_K+\mu $.
\end{proposition}
\begin{proof}
Determine $\ell$ so that $n/\rho^\ell\in T_u$. Then $\rho n/\rho^{\ell+1}\in T_u$ and 
\[\cpx{\rho n}_K=(\ell+1)\mu+d-u=\cpx{n}_K+\mu.\qedhere\]
\end{proof}

The above definition of $\cpx{m}_K$ may appear arbitrary. I have only two examples of arithmetic self-similar sets; in both cases this appears to be the correct election.

\begin{definition}\label{arithcompact}
An \emph{arithmetic compact} is a quasi-arithmetic self-similar compact set $K$ such that:
 
For any ordinal $0\le \alpha<\omega^\omega$ and $u\in\Z_{\ge0}$,  let $T_{u+1}[\alpha]=\frac{n}{\rho^k}$ with $\rho\nmid n$,  the two sets
\begin{gather}
\{T_{u}[\omega\alpha+r]\colon r\in\Z_{\ge0}\}, \\ \Bigl\{\frac{b(a \rho^r+1)}{\rho^{r+k}}\colon n=ab,\quad \cpx{n}_K=\cpx{a}_K+\cpx{b}_K,\quad r\in\Z_{\ge0}\Bigr\}\label{arcompfrac}
\end{gather}
have a finite symmetric difference. 
\end{definition}

\begin{remark}
The above definition is only provisional. The true definition will be that the quasi-arithmetic set will be arithmetic if we can describe the set $\{T_{u}[\omega\alpha+r]\colon r\in\Z_{\ge0}\}$ (at least except by a finite set) by some clear description. But definition \ref{arithcompact}, as given, only applies to one of the two examples that  I now have of arithmetical sets. 
\end{remark} 

\begin{proposition}\label{more}
Let $K$ be an arithmetic compact with module $\mu$ and ratio $\rho$. Then, for any natural number $n$ the $n$-th derivative set $K^{(n)}$ is also an arithmetic compact of the same module and ratio.
\end{proposition}

\begin{proof}
We only have to show that $K'$ is an arithmetical compact. Taking the derivative sets of  $\rho K^{(\mu)}=K$ we get $\rho(K')^{(\mu)}=K'$. $K'$ is compact because the set of limit points of a compact set is compact. Since it is contained in $K$, it is well ordered by $\preccurlyeq$. Since $K'=\{K[\omega\beta]\colon 1\le \beta\le\omega^\omega\}$, its order type is just $\omega^\omega+1$. So $K'$ is a self-similar compact set of ratio $\rho$ and module $\mu$. 

Our next goal is to prove that $\cpx{n}_K=\cpx{n}_{K'}$. Assume first that $1/\rho^a\in T_u=K^{(u)}\smallsetminus K^{(u+1)}=K'^{(u-1)}\smallsetminus K'^{(u)}$ with $0<u$ (there is always an $u>0$ satisfying this condition), then 
\[\cpx{1}_{K'}=a\mu+d'-(u-1)=0=\cpx{1}_K=a\mu+d-u.\]
It follows that the displacement $d'$ for $K'$ and that for $K$ are related by $d'=d-1$.

Now for any natural number $n$ not divisible by $\rho$ there is some $\ell$ and $u>0$ such that \[n/\rho^\ell\in K^{(u)}\smallsetminus K^{(u+1)}=K'^{(u-1)}\smallsetminus K'^{(u)}.\]
Therefore,
\[\cpx{n}_{K'}=\ell\mu+d'-(u-1)=\ell\mu+d-u=\cpx{n}_K.\]

Denote by $T_u(L)$ the set $T_u$ related to the self-similar compact set $L$. Then 
\[T_u(K')=(K')^{u}\smallsetminus (K')^{(u+1)}=T_{u+1}(K).\]
Therefore, the assertion of $K'$ being arithmetic is just part of the assertion of $K$ being arithmetic, thanks to the equality of the stable complexity associated to both compact sets.
\end{proof}

\section{The arithmetic compact set associated to the complexity of natural numbers}\label{existence}

Recall that we use $\cpx{n}$ and $\cpx{n}_\st$ as the complexity and the stable complexity of the natural number $n$. When $K$ is a quasi-aritmetic compact, $\cpx{n}_K$ is the corresponding stability gauge. We consider the sets defined in our conjectures \eqref{def:S}, but label them in a new way.
For $u\in\{0,1,2\}$, let 
\begin{equation}\label{Tcompl}
T_u:=\Bigl\{\frac{n}{3^k}\colon n \text{ stable and with } \cpx{n}=3k+2-u\Bigr\}.
\end{equation}
Hence $T_0=S_2$, $T_1=S_1$ and $T_2=S_0$. 
We will show that $K=\overline{T_0}$ is an arithmetic compact set with module $3$ and ratio $\rho=3$. Then we will have that these sets $T_u$ coincide with the sets $T_u(K)$ defined in Proposition \ref{P:cuatro}.  For the moment, in this section, consider only the sets $T_u$ defined in \eqref{Tcompl}.

\begin{theorem} The set $K=\overline{T_0}$ (with $T_0$ defined in \eqref{Tcompl}) is an arithmetic compact set with module $\mu=3$ and ratio $\rho=3$. The stability gauge $\cpx{n}_K$ coincides with the stable complexity $\cpx{n}_\st$.
\end{theorem}
\begin{proof}
The largest element with $\cpx{n}=3k+2$ is $E(k)=2\cdot 3^k$ \cite{intdft}*{Thm.~1.1}.
Hence if $n\in T_0$ with $\cpx{n}=3k+2$ we have $\frac{n}{3^k}\le 2$. Therefore, $T_0\subset[0,2]$ and the set $K$ is compact.

The set $T_u$ is related to the set $\mathscr{D}^{2-u}_\st$ in \cite{paperwo}*{Def.~2.4}
\[\mathscr{D}^{2-u}_\st:=\{\delta(n): n \text{ is a stable number with } \cpx{n}\equiv2-u\bmod3\}\]
With $f(x)=3^{(2-x)/3}$ we have $f(\mathscr{D}^{2-u}_\st+u)=T_u$. Then 
$f(\overline{\mathscr{D}^2_\st})$ is a closed set in $(0,\infty)$, with $0$ an accumulation point, because $\sup\mathscr{D}^2_\st=+\infty$. It follows that 
\[K\smallsetminus\{0\}=f(\overline{\mathscr{D}^2_\st})\]
The application $f$ is monotone decreasing. By \cite{paperwo}*{Thm.~7.4 (4)} the set $\overline{\mathscr{D}^2_\st}$ is well ordered, for the usual order, with order type $\omega^\omega$. Therefore, the set $K$  is well ordered by $\preccurlyeq$ with order type $\omega^\omega+1$, since it is  obtained by joining the image of this set with a limit element $0$ at the top. 

We have $(\overline{\mathscr{D}^2_\st})'''=\overline{\mathscr{D}^2_\st}+3$ (cfr. \cite{AA}*{Prop.~4.6}). By the above reasoning, this is equivalent to $(K\smallsetminus\{0\})'''=\frac13(K\smallsetminus\{0\})$. Adding  to this sets the element $\{0\}$, we get $K'''=\frac13K$. Therefore, $K$ is a self-similar compact with module $\mu=3$ and ratio $\rho=3$.

In \cite{AA}*{Cor.~4.7} it is shown that for $u\in\{0,1,2\}$
\begin{align*}
\overline{\mathscr{D}^2_\st}&=(\mathscr{D}^2_\st+3\Z_{\ge0})\cup(\mathscr{D}^1_\st+3\Z_{\ge0}+1)\cup(\mathscr{D}^0_\st+3\Z_{\ge0}+2)\\
\overline{\mathscr{D}^1_\st}&=(\mathscr{D}^1_\st+3\Z_{\ge0})\cup(\mathscr{D}^0_\st+3\Z_{\ge0}+1)\cup(\mathscr{D}^2_\st+3\Z_{\ge0}+2)\\
\overline{\mathscr{D}^0_\st}&=(\mathscr{D}^0_\st+3\Z_{\ge0})\cup(\mathscr{D}^2_\st+3\Z_{\ge0}+1)\cup(\mathscr{D}^1_\st+3\Z_{\ge0}+2).
\end{align*}
And by  \cite{AA}*{Prop.~4.6} we have 
\[(\overline{\mathscr{D}^2_\st})'=\overline{\mathscr{D}^1_\st}+1,\quad 
(\overline{\mathscr{D}^2_\st})''=\overline{\mathscr{D}^0_\st}+2.\]
Combining this
\begin{align*}
\overline{\mathscr{D}^2_\st}&=(\mathscr{D}^2_\st+3\Z_{\ge0})\cup(\mathscr{D}^1_\st+3\Z_{\ge0}+1)\cup(\mathscr{D}^0_\st+3\Z_{\ge0}+2)\\
(\overline{\mathscr{D}^2_\st})'&=(\mathscr{D}^1_\st+3\Z_{\ge0}+1)\cup(\mathscr{D}^0_\st+3\Z_{\ge0}+2)\cup(\mathscr{D}^2_\st+3\Z_{\ge0}+3)\\
(\overline{\mathscr{D}^2_\st})''&=(\mathscr{D}^0_\st+3\Z_{\ge0}+2)\cup(\mathscr{D}^2_\st+3\Z_{\ge0}+3)\cup(\mathscr{D}^1_\st+3\Z_{\ge0}+4).
\end{align*} 
Applying $f$ to these equations, we get the first three of the following equations 
\begin{align}\label{E:KT}
K\smallsetminus\{0\}&=\bigcup_{n\ge0}(T_0\cup T_1\cup T_2)/3^n\\
K'\smallsetminus\{0\}&=T_1\cup T_2\cup\bigcup_{n\ge1}(T_0\cup T_1\cup T_2)/3^n\nonumber\\
K''\smallsetminus\{0\}&=T_2\cup\bigcup_{n\ge1}(T_0\cup T_1\cup T_2)/3^n\nonumber\\
K'''\smallsetminus\{0\}&=\bigcup_{n\ge1}(T_0\cup T_1\cup T_2)/3^n\nonumber
\end{align}
the last equation follows from $K'''\smallsetminus\{0\}=(K\smallsetminus\{0\})/3$ and the first.
It is easy to show that the sets $T_u/3^n$ are disjoint. From this and the above formulas, we easily  get that 
\[T_0(K)=T_0,\quad T_1(K)=T_1,\quad T_2(K)=T_2.\]

To see that $K$ is quasi-arithmetic, we need to show (by \eqref{E:KT}) that there is a function $\kappa\colon(\N\smallsetminus 3\N)\to\Z_{\ge0}$ such that
\begin{equation}\label{E:quasiarith}
A:=\bigcup_{n\ge0}(T_0\cup T_1\cup T_2)/3^n=\Bigl\{\frac{m}{3^k}\colon 3\nmid m\text{ and }k\ge\kappa(m)\Bigr\}:=B.
\end{equation}

Let us define $\kappa$. Let $m$ be a natural number with $3\nmid m$, then by \cite{paper1}*{Thm.~13} there exist $a\ge0$ such that $n=m 3^a$ is stable. If $\cpx{n}=3\ell+2-u$ with $u\in\{0,1,2\}$, then $\frac{m}{3^{\ell-a}}\in T_u$. In this case, we define $\kappa(m)=\ell-a$. It is easy to see that this value does not depend on the chosen value of $a$, making $m 3^a$ stable. If $\cpx{n}=3\ell+2-u$ we have $n<3^{\ell+1}$, or $1\le m<3^{\ell-a+1}$, hence $\kappa(m)=\ell-a\ge0$.  Now, we prove \eqref{E:quasiarith}. 

It is clear that for this value of $m$ and any $k\ge \kappa(m)$ we have $m/3^k$ on the left hand side of \eqref{E:quasiarith}, and $B\subset A$. On the other hand, if $\frac{n}{3^k}\in A$, then  $\frac{n}{3^k}\in T_u/3^p$. It follows that
$\frac{n}{3^{k-p}}\in T_u$ so that $\frac{n}{3^{k-p}}=\frac{t}{3^s}$ with $t$ stable and 
$\cpx{t}=3s+2-u$, we will have $t=m 3^{a}$ with $3\nmid m$ and $\kappa(m)=s-a$. It follows that
\[\frac{n}{3^k}=\frac{1}{3^p}\frac{t}{3^s}=\frac{m3^a}{3^p 3^s}=\frac{m}{3^{p+s-a}}=\frac{m}{3^{p+\kappa(m)}}\]
with $p+\kappa(m)\ge\kappa(m)$. Therefore, $\frac{n}{3^k}\in B$. 
With this we see that $K$ is a quasi-arithmetic compact. 

Our next objective is to prove that $\cpx{n}_K=\cpx{n}_\st$. Since $3$ is a stable number with $\cpx{3}_\st=3$, we have $1\in T_2(K)=T_2=S_0$ and, therefore, 
$0=\cpx{1}_K=0\cdot3+d-2$ implies that the displacement $d=2$. For any $n$, there is always an $\ell$ and $u$ such that $n/3^\ell\in T_u(K)$ and in this case $\cpx{n}_K=3\ell+2-u$. If $n$ is stable $\cpx{n}=3\ell+2-u$. Hence, for $n$ stable $\cpx{n}_K=\cpx{n}_\st$.
When $n$ is not stable, there is $k$ with $n3^k$ stable, so that 
\[\cpx{n}_K+3k=\cpx{n 3^k}_K=\cpx{n3^k}_\st=\cpx{n}_\st+3k.\]

Showing that $K$ is arithmetic, given the equality of the gauge of $K$ and the stable complexity, is to show that the two sets
\begin{gather}
\{T_{u}[\omega\beta+r]\colon r\in\Z_{\ge0}\},\label{eq:8} \\ \Bigl\{\frac{b(a 3^r+1)}{3^{r+k}}\colon n=ab,\quad \cpx{n}_\st=\cpx{a}_\st+\cpx{b}_\st,\quad r\in\Z_{\ge0}\Bigr\}\label{eq:9}
\end{gather}
where $0\le \beta<\omega^\omega$,  $u\in\Z_{\ge0}$ and $T_{u+1}[\beta]=\frac{n}{3^k}$ with $3\nmid n$, has a finite symmetric difference.

We only have to show it for $u=0$, $1$ or $2$, because the  other cases differ only by dividing by a power of $3$. Take, for example, the case $u=2$. The hypothesis $T_3[\beta]=n/3^k$ is equivalent to $T_0[\beta]=n/3^{k-1}$, because $T_3=T_0/3$. This is equivalent to $S_2[\beta]=n/3^{k-1}$ so that $\cpx{n}_\st=3(k-1)+2=3k-1$, and $\dft_\st(n)=\mathscr{D}^2[\alpha]$. Proposition 4.19 in \cite{AA} implies  that the symmetric difference of the two sets
\[A:=\{\mathscr{D}^0_\st[\omega\beta+r]\colon r\ge0\},\quad B:= 
\{\dft_\st(b(a3^r+1)):r\ge0, n=ab, \cpx{n}_\st=\cpx{a}_\st+\cpx{b}_\st\}\]
is finite. 

The mapping $g(x)= 3^{-x/3} $ transform the set $A$ into $g(A)=\{S_0[\omega\beta+r]:r\ge0\}$. By Corollary 5.2 in \cite{AA} the numbers $b(a3^r+1)$ except for a finite sets satisfy 
\[\cpx{b(a3^r+1)}_\st=\cpx{a}_\st+\cpx{b}_\st+3r+1=3k+3r\]
Hence, the numbers in $B$ are of the form $\dft_\st(b(a3^r+1))=3k+3r-3\log_3(b(a3^r+1))$. It follows that the above map transforms the set $B$ into \eqref{eq:9} except for a finite set.  And we have proved that the symmetric difference between the sets \eqref{eq:8} and \eqref{eq:9} is finite in the case $u=2$. The other two cases are similar.
\end{proof}

\section{The unicity question}\label{unicity}

We have seen there is an arithmetic compact set $K$ of ratio and module equal 3. By Proposition \ref{more} $K^{(n)}$ are also arithmetic compact sets. Essentially there are three $K$, $K'$ and $K''$ and the others are homothetical to these ones. They have a distinct first element. So we may ask if there is a unique arithmetic compact set of ratio and module $3$ with $K[0]=2$. 

For each ordinal $\beta<\omega^\omega$ consider $T_1[\beta]=K[\omega(\beta+1)]$ this is the limit point of the sequence $K[\omega\beta+k]$, we have $K[\omega(\beta+1)]>K[\omega\beta]$ and we will call $I_\beta$ the interval in $\R$ $I_\beta:=(K[\omega(\beta+1)],K[\omega\beta])$ except when $\beta=0$ in which case we put 
$I_0=(K[\omega],K[0]]$. In this way, all elements in $T_0=K\smallsetminus K'$ are contained in some of these intervals. The intersection $I_\beta\cap T_0$ are called the $K$ sections. This section is a sequence with limit point $T_1[\beta]$. 

$K$ is an arithmetic compact, therefore, if $T_1[\beta]=\frac{n}{3^k}$, the terms of the section $I_\beta\cap T_0$ are given by \eqref{arcompfrac}.

The definition of the $I_\beta$ and sections can be extended to $K'$ and $K''$ in relation to $T_1$ and $T_2$.  The corresponding intervals will be called $I_\beta^u$. 

\begin{definition}\label{defsporadics}
Let $T_{u+1}[\beta]=\frac{n}{3^k}$. 
The numbers in the section $I_\beta^u\cap T_u$ not contained in the set 
\[A_\beta:=I_\beta^u\cap\Bigl\{\frac{b(a 3^r+1)}{3^{r+k}}\colon n=ab,\quad \cpx{n}_{\st}=\cpx{a}_\st+\cpx{b}_\st, \quad r\in\Z_{\ge0},\quad a^2+r^2>1\Bigr\},\]
will be called \emph{sporadics}.
\end{definition}

In the tables at the end of this paper I have written in boldface the fractions corresponding to sporadics.
For all sections given in the tables, we have $A_\beta\subset I_\beta\cap T_u$. This is not ensured by Theorem \ref{limitpoints}. I think this is a good definition of sporadics, but there are other possible ones leading to slightly different concepts\footnote{For example, in our tables at the end of this paper I arrange the fractions
taking $\frac{b(a 3^r+1)}{3^{r+k}}=\frac{n 3^{r+v}+b 3^v}{3^{r+k+v}}$, taking $v$ the greatest value with $b3^v < n$, and ordering the fractions according to decreasing values of $b 3^v$. This is the natural order of the fractions. In this case we could consider the values of $\frac{n 3^{\ell}+b 3^v}{3^{\ell+k}}$ for all $\ell\ge0$.\label{foot}}.

At the time of writing \cite{Arias}, uniqueness seemed very likely. The sporadics appear to have some structure, as very frequently they are powers of two. Theorem \ref{limitpoints} and determination of the sporadics, forced us the construction of the sets $T_u$. 
Conjectures 9, 10 and 11 were asking for the determination of these sporadics numbers. If these numbers followed a rule, I believed that the compact could be constructed in a unique way. In a sense, these three conjectures are not resolved. 

Hence, we have two questions. It is true that $A_\beta\subset I_\beta\cap T_u$?  Can we give a rule to determine the sporadics? Both questions appear to have a negative answer, but we do not know.

Regarding the rule, we have the following remark.

\subsection{Sporadics} By definition \ref{defsporadics}  the set $A_\beta\subset I_\beta$. The sporadics appear to be generated by the values of the fractions $a(b3^r+1)/3^{r+k}$ to the left of this interval. For example, take the section $K\cap(76/81,26/27)=(T_0[\omega^2+\omega2,T_0[\omega^2+\omega3])$ that contains many sporadics. Most of them appear from the sequences with limit $2048/2187$. Of the 19 sporadics in this section, 17 have this origin:
\begin{gather*}
\frac{18944}{19683}=\frac{512(4\cdot3^2+1)}{3^{2+7}},\quad \frac{6272}{6561}=\frac{128(16\cdot3+1)}{3^{1
+7}},\quad\frac{56320}{59049}=\frac{1024(2\cdot3^3+1)}{3^{3+7}},\\ \frac{2080}{2187}=\frac{32(64\cdot3^0+1)}{3^{0+7}},\quad \frac{18688}{19683}=\frac{256(8\cdot3^2+1)}{3^{2+7}},\quad
\frac{167936}{177147}=\frac{2048(3^4+1)}{3^{4+7}}\dots\\
\end{gather*}
The only ones that remain to be explained are the sporadic $\frac{700}{729}=\frac{25(3^3+1)}{3^{3+3}}$ and $\frac{2072}{2187}=\frac{56(4\cdot3^2+1)}{3^{2+5}}$.

This explains why the powers of two are usually sporadics. If $\cpx{2^n}=2n=3k+2-u$, then $\frac{2^n}{3^k}\in T_u$ and we have 
\[\frac{2^n}{3^k}=\frac{2^{n-2}(1\cdot 3+1)}{3^{(u-1)+1}},\quad\text{with}\quad \frac{2^{n-2}}{3^{u-1}}\in T_{u+1}.\]
Since the defects of the powers of two are relatively small, this makes $\frac{2^n}{3^k}$ usually  sporadic.

\subsection{H.~Altman almost counterexample to the unicity} H.~Altman thinks the unicity is not true.  To see it, he construct another example $H$ of an arithmetic compact with ratio and module 3. Although it seems very unlikely that $H=K$, we cannot prove $H\ne K$. We present here the example and then consider why $H=K$ appears unlikely.

He considers a modified complexity, instead of asking 
\[\cpx{ab}\le \cpx{a}+\cpx{b}, \quad  \cpx{a+b}\le \cpx{a}+\cpx{b}\]
his modified complexity $\cpx{\cdot}_1$ only satisfies
\[\cpx{ab}_1\le \cpx{a}_1+\cpx{b}_1, \quad  \cpx{a+1}_1\le \cpx{a}_1+1\]
The hole theory goes through, with $\cpx{n}_1$ in place of $\cpx{n}$ everywhere. There are just fewer low-defect polynomials (the ones with $+1$s only).  This gives us the new compact $H$. 

$H=K$ if and only if the stable complexities are equal $\cpx{n}_\st=\cpx{n}_{1,\st}$.

A solid number $n$ (see \cite{paper1} and \href{http://oeis.org/A195101}{OEIS:A195101})  is a number such that $n=a+b$ implies $\cpx{n}<\cpx{a}+\cpx{b}$. There are infinitely many solid numbers, the first are $1$, $6$, $9$, \dots And there are numbers $n$ whose complexity can only be computed with the decomposition $n=(n-b)+b$ where $b$ is a solid number. Hypothesis $K=H$ implies that these numbers that need a solid number $b>1$ are not stable. Usually, numbers of this form can be rewritten more efficiently. For example, numbers of the form 
$(4\cdot 3^k+1)+6$ do not have complexity $3k+10$ because 
\[(4\cdot 3^k+1)+6=6(2\cdot3^{k-1}+1)+1,\qquad k\ge1,\]
giving a complexity $\le 3k+6$. The first example in which we do not know how to rewrite it is $73(3^k+1)+6$, since $\cpx{73}=13$ we expect $\cpx{73(3^k+1)+6}=3k+19$. In fact, in \cite{deg1} we prove that for $k\ge k_0$ this is true. 

If $\cpx{73(3^k+1)+6}_\st\le 3k+ 18$, we will have for $\ell$ large 
\begin{align*}
\delta((73\cdot 3^k+79)3^\ell)&\le 3k+3\ell+18-3\log_3(73\cdot 3^k+79)3^\ell)\\ &=
18+3\log_3\frac{3^{k+\ell}}{73\cdot 3^{k+\ell}+79\cdot 3^\ell}
<18+3\log_3\frac{1}{73}\\ &\le 
6.284<58.614\cdot\delta(2)
\end{align*}
If we compute a covering set of $B_{59\delta(2)}$ we will know if this is true or not.
With the H. Altman algorithm, we have computed a good covering for $B_{42\delta(2)}$, but our Python programs crash every time I have tried to compute $B_{43\delta(2)}$. 

For $k=20$ the number  $\cpx{73(3^{20}+1)+6}=79$, but we obtain another representation
\[73(3^{20}+1)+6=2^3((((3^8+1)(3^6+1)3^2+1)(3^4+1)+1)3^2+1+1)\]
not using any solid number but 1, and giving its complexity.  

Using Iraids computations \cite{Ir}, for $k=21$ we have $\cpx{73(3^{21}+1)}=82$. This number is 
\[73(3^{21}+1)+6=763605783898=2\cdot 381802891949,\text{ with } \cpx{2}+\cpx{381802891949}=83\]
with $381802891949$ prime and $\cpx{763605783898-1}+\cpx{1}=23$. Therefore, in this case the complexity cannot be obtained by any other decomposition not using $b=6$. But what we do not know is that this number is stable. Its stability can only be proved by first computing a good covering of $B_{59\delta(2)}$. So we do not know.

\section{About the Tables of \texorpdfstring{$K$}{K}}
After the references, I add some tables on the compact $K$. They are constructed from the good covering of $B_{24\dft(2)}$.  With it we can construct all pairs $[m,3^{\lfloor\cpx{m}_\st\rfloor/3}]$ with $3\nmid m$ such that $\dft_\st(n)<24\dft(2)$. This means, for $\cpx{m}_\st=3k$, that $m/3^k>2^{24}/3^{16}$. So, all fractions $x$ in $T_u$ such that $x>3^{(2-u)/3}2^{24}/3^{16}$. In each case, we determine the sporadics and print them in boldface.  We terminate each section with the fractions in $3^k$ generating all non-sporadic terms ordered as explained in the footnote of page \pageref{foot}. In each case, the last numerical fractions given are the values of the fractions in $3^k$  for a given value of $k$. This will make it easy to prolong each section if needed.

\newpage

\section{Some elements of the compact \texorpdfstring{$K$: $T_0$}{K: T0}}
\tiny

{\setstretch{1.5}

\noindent $K[0]=2$,\hspace{4cm} $K\cap(\frac{4}{3},2]=[K[0],K[\omega])=[T_0[0], T_0[\omega])$\\
\noindent  $2$, $\frac{16}{9}$, $\frac{5}{3}$, $\mathbf{\frac{128}{81}}$, $\frac{14}{9}$, $\frac{40}{27}$, $\frac{13}{9}$, $\frac{38}{27}$, $\mathbf{\frac{1024}{729}}$, $\frac{112}{81}$, $\frac{37}{27}$, $\frac{110}{81}$, $\frac{328}{243}$, $\frac{109}{81}$, $\frac{326}{243}$, $\frac{976}{729}$, $\frac{325}{243}$, $\frac{974}{729}$, $\frac{2920}{2187}$, $\frac{973}{729}$, $\frac{2918}{2187}$, $\frac{8752}{6561}$, $\frac{2917}{2187}$, $\frac{8750}{6561}$, $\frac{26248}{19683}$, $\frac{8749}{6561}$, $\frac{26246}{19683}$, $\frac{78736}{59049}$, $\frac{26245}{19683}$, $\frac{78734}{59049}$, $\frac{236200}{177147}$, $\frac{78733}{59049}$, $\frac{236198}{177147}$, \dots$\frac{4\cdot3^k+4}{3^{k+1}}$, $\frac{4\cdot3^k+3}{3^{k+1}}$, $\frac{4\cdot3^k+2}{3^{k+1}}$, \dots

\noindent $K[\omega]=\frac{4}{3}$,\hspace{4cm} $K\cap(\frac{32}{27},\frac{4}{3})=[K[\omega+1],K[\omega2])=[T_0[\omega], T_0[\omega2])$\\
\noindent  $\frac{320}{243}$, $\mathbf{\frac{35}{27}}$, $\frac{104}{81}$, $\frac{34}{27}$, $\frac{304}{243}$, $\mathbf{\frac{8192}{6561}}$, $\frac{100}{81}$, $\frac{896}{729}$, $\frac{11}{9}$, $\frac{296}{243}$, $\frac{98}{81}$, $\frac{880}{729}$, $\frac{292}{243}$, $\frac{2624}{2187}$, $\frac{97}{81}$, $\frac{872}{729}$, $\frac{290}{243}$, $\frac{2608}{2187}$, $\frac{868}{729}$, $\frac{7808}{6561}$, $\frac{289}{243}$, $\frac{2600}{2187}$, $\frac{866}{729}$, $\frac{7792}{6561}$, $\frac{2596}{2187}$, $\frac{23360}{19683}$, $\frac{865}{729}$, $\frac{7784}{6561}$, $\frac{2594}{2187}$, $\frac{23344}{19683}$, $\frac{7780}{6561}$, \dots$\frac{32\cdot3^k+32}{3^{k+3}}$, $\frac{32\cdot3^k+27}{3^{k+3}}$, $\frac{32\cdot3^k+24}{3^{k+3}}$, $\frac{32\cdot3^k+18}{3^{k+3}}$, $\frac{32\cdot3^k+16}{3^{k+3}}$, $\frac{32\cdot3^k+12}{3^{k+3}}$, \dots

\noindent $K[\omega2]=\frac{32}{27}$,\hspace{4cm} $K\cap(\frac{10}{9},\frac{32}{27})=[K[\omega2+1],K[\omega3])=[T_0[\omega2], T_0[\omega3])$\\
\noindent  $\frac{95}{81}$, $\mathbf{\frac{2560}{2187}}$, $\frac{280}{243}$, $\frac{31}{27}$, $\mathbf{\frac{832}{729}}$, $\frac{92}{81}$, $\frac{275}{243}$, $\frac{820}{729}$, $\frac{91}{81}$, $\frac{272}{243}$, $\frac{815}{729}$, $\frac{2440}{2187}$, $\frac{271}{243}$, $\frac{812}{729}$, $\frac{2435}{2187}$, $\frac{7300}{6561}$, $\frac{811}{729}$, $\frac{2432}{2187}$, $\frac{7295}{6561}$, $\frac{21880}{19683}$, $\frac{2431}{2187}$, $\frac{7292}{6561}$, $\frac{21875}{19683}$, $\frac{65620}{59049}$, $\frac{7291}{6561}$, $\frac{21872}{19683}$, $\frac{65615}{59049}$, $\frac{196840}{177147}$, $\frac{21871}{19683}$, $\frac{65612}{59049}$, $\frac{196835}{177147}$, \dots$\frac{10\cdot3^k+10}{3^{k+2}}$, $\frac{10\cdot3^k+9}{3^{k+2}}$, $\frac{10\cdot3^k+6}{3^{k+2}}$, $\frac{10\cdot3^k+5}{3^{k+2}}$, \dots

\noindent $K[\omega3]=\frac{10}{9}$,\hspace{4cm} $K\cap(\frac{256}{243},\frac{10}{9})=[K[\omega3+1],K[\omega4])=[T_0[\omega3], T_0[\omega4])$\\
\noindent  $\mathbf{\frac{65536}{59049}}$, $\frac{800}{729}$, $\mathbf{\frac{266}{243}}$, $\frac{7168}{6561}$, $\frac{88}{81}$, $\frac{2368}{2187}$, $\frac{784}{729}$, $\mathbf{\frac{29}{27}}$, $\frac{7040}{6561}$, $\frac{260}{243}$, $\frac{2336}{2187}$, $\frac{20992}{19683}$, $\mathbf{\frac{259}{243}}$, $\frac{776}{729}$, $\frac{6976}{6561}$, $\frac{86}{81}$, $\frac{2320}{2187}$, $\frac{20864}{19683}$, $\frac{772}{729}$, $\frac{6944}{6561}$, $\frac{62464}{59049}$, $\frac{257}{243}$, $\frac{2312}{2187}$, $\frac{20800}{19683}$, $\frac{770}{729}$, $\frac{6928}{6561}$, $\frac{62336}{59049}$, $\frac{2308}{2187}$, $\frac{20768}{19683}$, $\frac{186880}{177147}$, $\frac{769}{729}$, $\frac{6920}{6561}$, $\frac{62272}{59049}$, $\frac{2306}{2187}$, $\frac{20752}{19683}$, $\frac{186752}{177147}$, $\frac{6916}{6561}$, $\frac{62240}{59049}$, \dots$\frac{256\cdot3^k+256}{3^{k+5}}$, $\frac{256\cdot3^k+243}{3^{k+5}}$, $\frac{256\cdot3^k+216}{3^{k+5}}$, $\frac{256\cdot3^k+192}{3^{k+5}}$, $\frac{256\cdot3^k+162}{3^{k+5}}$, $\frac{256\cdot3^k+144}{3^{k+5}}$, $\frac{256\cdot3^k+128}{3^{k+5}}$, $\frac{256\cdot3^k+108}{3^{k+5}}$, $\frac{256\cdot3^k+96}{3^{k+5}}$, \dots

\noindent $K[\omega4]=\frac{256}{243}$,\hspace{4cm} $K\cap(\frac{28}{27},\frac{256}{243})=[K[\omega4+1],K[\omega5])=[T_0[\omega4], T_0[\omega5])$\\
\noindent  $\frac{2296}{2187}$, $\frac{85}{81}$, $\frac{763}{729}$, $\frac{254}{243}$, $\frac{2282}{2187}$, $\frac{760}{729}$, $\frac{6832}{6561}$, $\frac{253}{243}$, $\mathbf{\frac{20480}{19683}}$, $\frac{2275}{2187}$, $\frac{758}{729}$, $\frac{6818}{6561}$, $\frac{2272}{2187}$, $\frac{20440}{19683}$, $\frac{757}{729}$, $\frac{6811}{6561}$, $\frac{2270}{2187}$, $\frac{20426}{19683}$, $\frac{6808}{6561}$, $\frac{61264}{59049}$, $\frac{2269}{2187}$, $\frac{20419}{19683}$, $\frac{6806}{6561}$, $\frac{61250}{59049}$, $\frac{20416}{19683}$,  \dots$\frac{28\cdot3^k+28}{3^{k+3}}$, $\frac{28\cdot3^k+27}{3^{k+3}}$, $\frac{28\cdot3^k+21}{3^{k+3}}$, $\frac{28\cdot3^k+18}{3^{k+3}}$, $\frac{28\cdot3^k+14}{3^{k+3}}$, $\frac{28\cdot3^k+12}{3^{k+3}}$, \dots

\noindent $K[\omega5]=\frac{28}{27}$,\hspace{4cm} $K\cap(\frac{82}{81},\frac{28}{27})=[K[\omega5+1],K[\omega6])=[T_0[\omega5], T_0[\omega6])$\\
\noindent  $\mathbf{\frac{250}{243}}$, $\frac{6724}{6561}$, $\frac{83}{81}$, $\mathbf{\frac{2240}{2187}}$, $\mathbf{\frac{248}{243}}$, $\frac{20008}{19683}$, $\frac{247}{243}$, $\mathbf{\frac{740}{729}}$, $\mathbf{\frac{6656}{6561}}$, $\frac{59860}{59049}$, $\frac{739}{729}$, $\frac{179416}{177147}$, $\frac{2215}{2187}$, $\frac{538084}{531441}$, $\frac{6643}{6561}$, $\frac{1614088}{1594323}$, $\frac{19927}{19683}$, $\frac{4842100}{4782969}$, $\frac{59779}{59049}$, $\frac{14526136}{14348907}$, $\frac{179335}{177147}$, $\frac{43578244}{43046721}$, $\frac{538003}{531441}$, $\frac{130734568}{129140163}$, $\frac{1614007}{1594323}$, \dots$\frac{82\cdot3^k+82}{3^{k+4}}$, $\frac{82\cdot3^k+81}{3^{k+4}}$, \dots

\noindent $K[\omega6]=\frac{82}{81}$,\hspace{4cm} $K\cap(\frac{244}{243},\frac{82}{81})=[K[\omega6+1],K[\omega7])=[T_0[\omega6], T_0[\omega7])$\\
\noindent  $\mathbf{\frac{736}{729}}$, $\frac{59536}{59049}$, $\frac{245}{243}$, $\mathbf{\frac{2200}{2187}}$, $\frac{178120}{177147}$, $\frac{733}{729}$, $\frac{533872}{531441}$, $\frac{2197}{2187}$, $\frac{1601128}{1594323}$, $\frac{6589}{6561}$, $\frac{4802896}{4782969}$, $\frac{19765}{19683}$, $\frac{14408200}{14348907}$, $\frac{59293}{59049}$, $\frac{43224112}{43046721}$, $\frac{177877}{177147}$, $\frac{129671848}{129140163}$, $\frac{533629}{531441}$, \dots$\frac{244\cdot3^k+244}{3^{k+5}}$, $\frac{244\cdot3^k+243}{3^{k+5}}$, \dots

\noindent $K[\omega7]=\frac{244}{243}$,\hspace{4cm} $K\cap(\frac{730}{729},\frac{244}{243})=[K[\omega7+1],K[\omega8])=[T_0[\omega7], T_0[\omega8])$\\
\noindent  $\frac{532900}{531441}$, $\frac{731}{729}$, $\frac{1597240}{1594323}$, $\frac{2191}{2187}$, $\frac{4790260}{4782969}$, $\frac{6571}{6561}$, $\frac{14369320}{14348907}$, $\frac{19711}{19683}$, $\frac{43106500}{43046721}$, $\frac{59131}{59049}$, $\frac{129318040}{129140163}$, $\frac{177391}{177147}$,\dots$\frac{730\cdot3^k+730}{3^{k+6}}$, $\frac{730\cdot3^k+729}{3^{k+6}}$, \dots

\dots

\noindent $K[\omega^2]=1$,\hspace{4cm} $K\cap(\frac{80}{81},1)=[K[\omega^2+1],K[\omega^2+\omega])=[T_0[\omega^2], T_0[\omega^2+\omega])$\\
\noindent  $\frac{6560}{6561}$, $\frac{728}{729}$, $\frac{2180}{2187}$, $\frac{242}{243}$, $\frac{2176}{2187}$, $\frac{725}{729}$, $\frac{6520}{6561}$, $\frac{724}{729}$, $\frac{2170}{2187}$, $\frac{241}{243}$, $\frac{19520}{19683}$, $\frac{2168}{2187}$, $\frac{6500}{6561}$, $\frac{722}{729}$, $\frac{6496}{6561}$, $\frac{2165}{2187}$, $\frac{19480}{19683}$, $\frac{2164}{2187}$, $\frac{6490}{6561}$, $\frac{721}{729}$, $\frac{58400}{59049}$, $\frac{6488}{6561}$, $\frac{19460}{19683}$, $\frac{2162}{2187}$, $\frac{19456}{19683}$, $\frac{6485}{6561}$, $\frac{58360}{59049}$, $\frac{6484}{6561}$, $\frac{19450}{19683}$, $\frac{2161}{2187}$,\dots$\frac{80\cdot3^k+80}{3^{k+4}}$, $\frac{80\cdot3^k+72}{3^{k+4}}$, $\frac{80\cdot3^k+60}{3^{k+4}}$, $\frac{80\cdot3^k+54}{3^{k+4}}$, $\frac{80\cdot3^k+48}{3^{k+4}}$, $\frac{80\cdot3^k+45}{3^{k+4}}$, $\frac{80\cdot3^k+40}{3^{k+4}}$, $\frac{80\cdot3^k+36}{3^{k+4}}$, $\frac{80\cdot3^k+30}{3^{k+4}}$, $\frac{80\cdot3^k+27}{3^{k+4}}$, \dots

\noindent $K[\omega^2+\omega]=\frac{80}{81}$,\hspace{4cm} $K\cap(\frac{26}{27},\frac{80}{81})=[K[\omega^2+\omega+1],K[\omega^2+\omega2])=[T_0[\omega^2+\omega], T_0[\omega^2+\omega2])$\\
\noindent  $\mathbf{\frac{524288}{531441}}$, $\frac{715}{729}$, $\mathbf{\frac{238}{243}}$, $\mathbf{\frac{6400}{6561}}$, $\frac{79}{81}$, $\frac{2132}{2187}$, $\mathbf{\frac{2128}{2187}}$, $\frac{236}{243}$, $\mathbf{\frac{57344}{59049}}$, $\frac{2119}{2187}$, $\frac{235}{243}$, $\frac{6344}{6561}$, $\frac{704}{729}$, $\frac{6331}{6561}$, $\frac{703}{729}$, $\frac{18980}{19683}$, $\frac{2108}{2187}$, $\frac{18967}{19683}$, $\frac{2107}{2187}$, $\frac{56888}{59049}$, $\frac{6320}{6561}$, $\frac{56875}{59049}$, $\frac{6319}{6561}$, $\frac{170612}{177147}$, $\frac{18956}{19683}$, $\frac{170599}{177147}$, $\frac{18955}{19683}$, $\frac{511784}{531441}$, $\frac{56864}{59049}$, $\frac{511771}{531441}$, $\frac{56863}{59049}$, \dots$\frac{26\cdot3^k+26}{3^{k+3}}$, $\frac{26\cdot3^k+18}{3^{k+3}}$, $\frac{26\cdot3^k+13}{3^{k+3}}$, $\frac{26\cdot3^k+9}{3^{k+3}}$, \dots

\noindent $K[\omega^2+\omega2]=\frac{26}{27}$,\hspace{4cm} $K\cap(\frac{76}{81},\frac{26}{27})=[K[\omega^2+\omega2+1],K[\omega^2+\omega3])=[T_0[\omega^2+\omega2], T_0[\omega^2+\omega3])$\\
\noindent  $\mathbf{\frac{18944}{19683}}$, $\mathbf{\frac{700}{729}}$, $\mathbf{\frac{6272}{6561}}$, $\frac{2090}{2187}$, $\frac{232}{243}$, $\mathbf{\frac{56320}{59049}}$, $\mathbf{\frac{2080}{2187}}$, $\frac{77}{81}$, $\frac{6232}{6561}$, $\mathbf{\frac{18688}{19683}}$, $\mathbf{\frac{167936}{177147}}$, $\mathbf{\frac{2072}{2187}}$, $\frac{2071}{2187}$, $\frac{230}{243}$, $\mathbf{\frac{6208}{6561}}$, $\mathbf{\frac{55808}{59049}}$, $\frac{6194}{6561}$, $\frac{688}{729}$, $\mathbf{\frac{18560}{19683}}$, $\frac{229}{243}$, $\mathbf{\frac{166912}{177147}}$, $\frac{18544}{19683}$, $\mathbf{\frac{6176}{6561}}$, $\frac{6175}{6561}$, $\frac{686}{729}$, $\mathbf{\frac{55552}{59049}}$, $\mathbf{\frac{499712}{531441}}$, $\frac{18506}{19683}$, $\frac{2056}{2187}$, $\mathbf{\frac{18496}{19683}}$, $\frac{685}{729}$, $\frac{55480}{59049}$, $\mathbf{\frac{166400}{177147}}$, $\frac{18487}{19683}$, $\frac{2054}{2187}$, $\frac{55442}{59049}$, $\frac{6160}{6561}$, $\frac{2053}{2187}$, $\frac{166288}{177147}$, $\mathbf{\frac{55424}{59049}}$, $\frac{55423}{59049}$, $\frac{6158}{6561}$, $\frac{166250}{177147}$, $\frac{18472}{19683}$, $\frac{6157}{6561}$, $\frac{498712}{531441}$, $\frac{166231}{177147}$, $\frac{18470}{19683}$, $\mathbf{\frac{498688}{531441}}$, $\frac{498674}{531441}$, $\frac{55408}{59049}$, $\frac{18469}{19683}$, $\frac{1495984}{1594323}$, $\frac{498655}{531441}$, $\frac{55406}{59049}$, $\frac{1495946}{1594323}$, $\frac{166216}{177147}$, $\frac{55405}{59049}$, $\frac{4487800}{4782969}$, $\frac{1495927}{1594323}$, $\frac{166214}{177147}$, $\frac{4487762}{4782969}$, $\frac{498640}{531441}$, $\frac{166213}{177147}$, $\frac{13463248}{14348907}$, $\frac{4487743}{4782969}$, $\frac{498638}{531441}$, $\frac{13463210}{14348907}$, $\frac{1495912}{1594323}$, $\frac{498637}{531441}$, \dots$\frac{76\cdot3^k+76}{3^{k+4}}$, $\frac{76\cdot3^k+57}{3^{k+4}}$, $\frac{76\cdot3^k+54}{3^{k+4}}$, $\frac{76\cdot3^k+38}{3^{k+4}}$, $\frac{76\cdot3^k+36}{3^{k+4}}$, $\frac{76\cdot3^k+27}{3^{k+4}}$, \dots

\noindent $K[\omega^2+\omega3]=\frac{76}{81}$,\hspace{4cm} $K\cap(\frac{2048}{2187},\frac{76}{81})=[K[\omega^2+\omega3+1],K[\omega^2+\omega4])=[T_0[\omega^2+\omega3], T_0[\omega^2+\omega4])$\\
\noindent  $\frac{18464}{19683}$, $\frac{166144}{177147}$, $\frac{1495040}{1594323}$, $\frac{6152}{6561}$, $\frac{55360}{59049}$, $\frac{498176}{531441}$, $\frac{2050}{2187}$, $\frac{18448}{19683}$, $\frac{166016}{177147}$, $\frac{1494016}{1594323}$, $\frac{6148}{6561}$, $\frac{55328}{59049}$, $\frac{497920}{531441}$, $\frac{683}{729}$, $\frac{4481024}{4782969}$, $\frac{18440}{19683}$, $\frac{165952}{177147}$, $\frac{1493504}{1594323}$, $\frac{6146}{6561}$, $\frac{55312}{59049}$, $\frac{497792}{531441}$, $\frac{4480000}{4782969}$, $\frac{18436}{19683}$, $\frac{165920}{177147}$, $\frac{1493248}{1594323}$, $\frac{6145}{6561}$, \dots$\frac{2048\cdot3^k+2048}{3^{k+7}}$, $\frac{2048\cdot3^k+1944}{3^{k+7}}$, $\frac{2048\cdot3^k+1728}{3^{k+7}}$, $\frac{2048\cdot3^k+1536}{3^{k+7}}$,
 $\frac{2048\cdot3^k+1458}{3^{k+7}}$, $\frac{2048\cdot3^k+1296}{3^{k+7}}$, $\frac{2048\cdot3^k+1152}{3^{k+7}}$, $\frac{2048\cdot3^k+1024}{3^{k+7}}$, $\frac{2048\cdot3^k+972}{3^{k+7}}$, $\frac{2048\cdot3^k+864}{3^{k+7}}$, $\frac{2048\cdot3^k+768}{3^{k+7}}$, $\frac{2048\cdot3^k+729}{3^{k+7}}$, \dots

\noindent $K[\omega^2+\omega4]=\frac{2048}{2187}$,\hspace{3.cm} $K\cap(\frac{25}{27},\frac{2048}{2187})=[K[\omega^2+\omega4+1],K[\omega^2+\omega5])=[T_0[\omega^2+\omega4], T_0[\omega^2+\omega5])$\\
\noindent  $\mathbf{\frac{2044}{2187}}$, $\mathbf{\frac{18368}{19683}}$, $\frac{680}{729}$, $\mathbf{\frac{679}{729}}$, $\mathbf{\frac{2035}{2187}}$, $\mathbf{\frac{6104}{6561}}$, $\frac{226}{243}$, $\frac{6100}{6561}$, $\mathbf{\frac{2032}{2187}}$, $\frac{2030}{2187}$, $\mathbf{\frac{18256}{19683}}$, $\frac{676}{729}$, $\frac{18250}{19683}$, $\frac{6080}{6561}$, $\frac{2026}{2187}$, $\frac{54700}{59049}$, $\frac{18230}{19683}$, $\frac{6076}{6561}$, $\frac{164050}{177147}$, $\frac{54680}{59049}$, $\frac{18226}{19683}$, $\frac{492100}{531441}$, $\frac{164030}{177147}$, $\frac{54676}{59049}$, $\frac{1476250}{1594323}$, $\frac{492080}{531441}$, $\frac{164026}{177147}$, $\frac{4428700}{4782969}$, $\frac{1476230}{1594323}$, $\frac{492076}{531441}$, \dots$\frac{25\cdot3^k+25}{3^{k+3}}$, $\frac{25\cdot3^k+15}{3^{k+3}}$, $\frac{25\cdot3^k+9}{3^{k+3}}$, \dots

\dots

\noindent $K[\omega^22]=\frac{8}{9}$,\hspace{4cm} $K\cap(\frac{640}{729},\frac{8}{9})=[K[\omega^22+1],K[\omega^22+\omega])=[T_0[\omega^2 2], T_0[\omega^2 2+\omega])$\\
\noindent  $\frac{52480}{59049}$, $\frac{5824}{6561}$, $\frac{1940}{2187}$, $\mathbf{\frac{646}{729}}$, $\frac{17440}{19683}$, $\frac{1936}{2187}$, $\frac{215}{243}$, $\frac{17408}{19683}$, $\frac{5800}{6561}$, $\frac{644}{729}$, $\frac{52160}{59049}$, $\frac{5792}{6561}$, $\frac{1930}{2187}$, $\frac{17360}{19683}$, $\frac{1928}{2187}$, $\frac{156160}{177147}$, $\frac{17344}{19683}$, $\frac{5780}{6561}$, $\frac{214}{243}$, $\frac{52000}{59049}$, $\frac{5776}{6561}$, $\frac{1925}{2187}$, $\frac{51968}{59049}$, $\frac{17320}{19683}$, $\frac{1924}{2187}$, $\frac{155840}{177147}$, $\frac{17312}{19683}$, $\frac{5770}{6561}$, $\frac{641}{729}$, $\frac{51920}{59049}$, $\frac{5768}{6561}$, \dots$\frac{640\cdot3^k+640}{3^{k+6}}$, $\frac{640\cdot3^k+576}{3^{k+6}}$, $\frac{640\cdot3^k+540}{3^{k+6}}$, $\frac{640\cdot3^k+486}{3^{k+6}}$, $\frac{640\cdot3^k+480}{3^{k+6}}$, $\frac{640\cdot3^k+432}{3^{k+6}}$, $\frac{640\cdot3^k+405}{3^{k+6}}$, $\frac{640\cdot3^k+384}{3^{k+6}}$, $\frac{640\cdot3^k+360}{3^{k+6}}$, $\frac{640\cdot3^k+324}{3^{k+6}}$, $\frac{640\cdot3^k+320}{3^{k+6}}$, $\frac{640\cdot3^k+288}{3^{k+6}}$, $\frac{640\cdot3^k+270}{3^{k+6}}$, $\frac{640\cdot3^k+243}{3^{k+6}}$, $\frac{640\cdot3^k+240}{3^{k+6}}$, $\frac{640\cdot3^k+216}{3^{k+6}}$, \dots

\noindent $K[\omega^22+\omega]=\frac{640}{729}$,\hspace{3.5cm} $K\cap(\frac{70}{81},\frac{640}{729})=[K[\omega^22+\omega+1],K[\omega^22+\omega2])=[T_0[\omega^2 2+\omega], T_0[\omega^2 2+\omega2])$\\
\noindent  $\mathbf{\frac{4194304}{4782969}}$, $\frac{71}{81}$, $\frac{5740}{6561}$, $\frac{637}{729}$, $\frac{212}{243}$, $\mathbf{\frac{5720}{6561}}$, $\frac{635}{729}$, $\frac{1904}{2187}$, $\frac{5705}{6561}$, $\frac{1900}{2187}$, $\frac{211}{243}$, $\mathbf{\frac{1898}{2187}}$, $\frac{17080}{19683}$, $\frac{1897}{2187}$, $\mathbf{\frac{51200}{59049}}$, $\frac{632}{729}$, $\mathbf{\frac{17056}{19683}}$, $\frac{1895}{2187}$, $\frac{5684}{6561}$, $\frac{17045}{19683}$, $\frac{5680}{6561}$, $\frac{631}{729}$, $\frac{51100}{59049}$, $\frac{5677}{6561}$, $\frac{1892}{2187}$, $\frac{5675}{6561}$, $\frac{17024}{19683}$, $\frac{51065}{59049}$, $\frac{17020}{19683}$, $\frac{1891}{2187}$, $\frac{153160}{177147}$, $\frac{17017}{19683}$, $\frac{5672}{6561}$, $\frac{17015}{19683}$, $\frac{51044}{59049}$, $\frac{153125}{177147}$, $\frac{51040}{59049}$, $\frac{5671}{6561}$, \dots$\frac{70\cdot3^k+70}{3^{k+4}}$, $\frac{70\cdot3^k+63}{3^{k+4}}$, $\frac{70\cdot3^k+54}{3^{k+4}}$, $\frac{70\cdot3^k+45}{3^{k+4}}$, $\frac{70\cdot3^k+42}{3^{k+4}}$, $\frac{70\cdot3^k+35}{3^{k+4}}$, $\frac{70\cdot3^k+30}{3^{k+4}}$, $\frac{70\cdot3^k+27}{3^{k+4}}$, \dots

\noindent $K[\omega^22+\omega2]=\frac{70}{81}$,\hspace{3.5cm} $K\cap(\frac{208}{243},\frac{70}{81})=[K[\omega^22+\omega2+1],K[\omega^22+\omega3])=[T_0[\omega^2 2+\omega2], T_0[\omega^2 2+\omega3])$\\
\noindent  $\frac{5668}{6561}$, $\frac{1888}{2187}$, $\mathbf{\frac{458752}{531441}}$, $\mathbf{\frac{629}{729}}$, $\frac{1885}{2187}$, $\frac{628}{729}$, $\frac{16952}{19683}$, $\frac{209}{243}$, $\frac{5642}{6561}$, $\frac{1880}{2187}$, $\frac{50752}{59049}$, $\frac{626}{729}$, $\frac{16900}{19683}$, $\frac{5632}{6561}$, $\frac{5629}{6561}$, $\frac{1876}{2187}$, $\frac{50648}{59049}$, $\frac{625}{729}$, $\frac{16874}{19683}$, $\frac{5624}{6561}$, $\frac{151840}{177147}$, $\frac{1874}{2187}$, $\frac{50596}{59049}$, $\frac{16864}{19683}$, $\frac{16861}{19683}$, $\frac{5620}{6561}$, $\frac{151736}{177147}$, $\frac{1873}{2187}$, $\frac{50570}{59049}$, $\frac{16856}{19683}$, \dots$\frac{208\cdot3^k+208}{3^{k+5}}$, $\frac{208\cdot3^k+162}{3^{k+5}}$, $\frac{208\cdot3^k+156}{3^{k+5}}$, $\frac{208\cdot3^k+144}{3^{k+5}}$, $\frac{208\cdot3^k+117}{3^{k+5}}$, $\frac{208\cdot3^k+108}{3^{k+5}}$, $\frac{208\cdot3^k+104}{3^{k+5}}$, $\frac{208\cdot3^k+81}{3^{k+5}}$, $\frac{208\cdot3^k+78}{3^{k+5}}$, $\frac{208\cdot3^k+72}{3^{k+5}}$, \dots

\noindent $K[\omega^22+\omega3]=\frac{208}{243}$,\hspace{3.5cm} $K\cap(\frac{68}{81},\frac{208}{243})=[K[\omega^22+\omega3+1],K[\omega^22+\omega4])=[T_0[\omega^2 2+\omega3], T_0[\omega^2 2+\omega4])$\\
\noindent  $\mathbf{\frac{151552}{177147}}$, $\frac{1870}{2187}$, $\mathbf{\frac{5600}{6561}}$, $\frac{23}{27}$, $\mathbf{\frac{1862}{2187}}$, $\mathbf{\frac{620}{729}}$, $\frac{5576}{6561}$, $\mathbf{\frac{50176}{59049}}$, $\mathbf{\frac{16720}{19683}}$, $\mathbf{\frac{1856}{2187}}$, $\mathbf{\frac{450560}{531441}}$, $\frac{206}{243}$, $\frac{1853}{2187}$, $\mathbf{\frac{1850}{2187}}$, $\mathbf{\frac{5548}{6561}}$, $\mathbf{\frac{16640}{19683}}$, $\frac{616}{729}$, $\frac{5542}{6561}$, $\mathbf{\frac{49856}{59049}}$, $\mathbf{\frac{149504}{177147}}$, $\frac{205}{243}$, $\frac{16592}{19683}$, $\mathbf{\frac{1843}{2187}}$, $\mathbf{\frac{1343488}{1594323}}$, $\frac{614}{729}$, $\mathbf{\frac{16576}{19683}}$, $\frac{5525}{6561}$, $\mathbf{\frac{16568}{19683}}$, $\frac{1840}{2187}$, $\frac{16558}{19683}$, $\mathbf{\frac{49664}{59049}}$, $\frac{613}{729}$, $\frac{49640}{59049}$, $\frac{1838}{2187}$, $\frac{16541}{19683}$, $\frac{5512}{6561}$, $\mathbf{\frac{446464}{531441}}$, $\frac{49606}{59049}$, $\frac{1837}{2187}$, $\frac{148784}{177147}$, $\frac{5510}{6561}$, $\frac{49589}{59049}$, $\frac{16528}{19683}$, $\frac{148750}{177147}$, $\frac{5509}{6561}$, $\frac{446216}{531441}$, $\frac{16526}{19683}$, $\frac{148733}{177147}$, $\frac{49576}{59049}$, $\frac{446182}{531441}$, $\frac{16525}{19683}$, $\frac{1338512}{1594323}$, $\frac{49574}{59049}$, $\frac{446165}{531441}$, $\frac{148720}{177147}$, $\frac{1338478}{1594323}$, $\frac{49573}{59049}$, \dots$\frac{68\cdot3^k+68}{3^{k+4}}$, $\frac{68\cdot3^k+54}{3^{k+4}}$, $\frac{68\cdot3^k+51}{3^{k+4}}$, $\frac{68\cdot3^k+36}{3^{k+4}}$, $\frac{68\cdot3^k+34}{3^{k+4}}$, $\frac{68\cdot3^k+27}{3^{k+4}}$, \dots

\noindent $K[\omega^22+\omega4]=\frac{68}{81}$,\hspace{3.5cm} $K\cap(\frac{608}{729},\frac{68}{81})=[K[\omega^22+\omega4+1],K[\omega^22+\omega5])=[T_0[\omega^2 2+\omega4], T_0[\omega^2 2+\omega5])$\\
\noindent  $\frac{49552}{59049}$, $\frac{5504}{6561}$, $\mathbf{\frac{5500}{6561}}$, $\mathbf{\frac{148480}{177147}}$, $\frac{16492}{19683}$, $\frac{1832}{2187}$, $\mathbf{\frac{1335296}{1594323}}$, $\frac{148352}{177147}$, $\frac{5491}{6561}$, $\frac{610}{729}$, $\mathbf{\frac{49408}{59049}}$, $\frac{49400}{59049}$, $\frac{5488}{6561}$, $\mathbf{\frac{444416}{531441}}$, $\frac{16454}{19683}$, $\frac{1828}{2187}$, $\mathbf{\frac{3997696}{4782969}}$, $\frac{148048}{177147}$, $\frac{16448}{19683}$, $\frac{203}{243}$, $\frac{49324}{59049}$, $\mathbf{\frac{147968}{177147}}$, $\frac{5480}{6561}$, $\frac{443840}{531441}$, $\frac{16435}{19683}$, $\mathbf{\frac{1331200}{1594323}}$, $\frac{1826}{2187}$, $\frac{147896}{177147}$, $\frac{16432}{19683}$, $\frac{49286}{59049}$, $\frac{5476}{6561}$, $\frac{443536}{531441}$, $\frac{49280}{59049}$, $\frac{1825}{2187}$, $\frac{147820}{177147}$, $\frac{16424}{19683}$, $\frac{1330304}{1594323}$, $\frac{49267}{59049}$, $\frac{5474}{6561}$, $\mathbf{\frac{443392}{531441}}$, $\frac{443384}{531441}$, $\frac{49264}{59049}$, $\frac{147782}{177147}$, $\frac{16420}{19683}$, $\frac{1330000}{1594323}$, $\frac{147776}{177147}$, $\frac{5473}{6561}$, $\frac{443308}{531441}$, $\frac{49256}{59049}$, $\frac{3989696}{4782969}$, $\frac{147763}{177147}$, $\frac{16418}{19683}$, $\frac{1329848}{1594323}$, $\frac{147760}{177147}$, $\mathbf{\frac{3989504}{4782969}}$, $\frac{443270}{531441}$, $\frac{49252}{59049}$, $\frac{3989392}{4782969}$, $\frac{443264}{531441}$, $\frac{16417}{19683}$, $\frac{1329772}{1594323}$, $\frac{147752}{177147}$, $\frac{11967872}{14348907}$, $\frac{443251}{531441}$, $\frac{49250}{59049}$, $\frac{3989240}{4782969}$, $\frac{443248}{531441}$, $\frac{1329734}{1594323}$, $\frac{147748}{177147}$, $\frac{11967568}{14348907}$, $\frac{1329728}{1594323}$, $\frac{49249}{59049}$, $\frac{3989164}{4782969}$, $\frac{443240}{531441}$, \dots$\frac{608\cdot3^k+608}{3^{k+6}}$, $\frac{608\cdot3^k+513}{3^{k+6}}$, $\frac{608\cdot3^k+486}{3^{k+6}}$, $\frac{608\cdot3^k+456}{3^{k+6}}$, $\frac{608\cdot3^k+432}{3^{k+6}}$, $\frac{608\cdot3^k+342}{3^{k+6}}$, $\frac{608\cdot3^k+324}{3^{k+6}}$, $\frac{608\cdot3^k+304}{3^{k+6}}$, $\frac{608\cdot3^k+288}{3^{k+6}}$, $\frac{608\cdot3^k+243}{3^{k+6}}$, $\frac{608\cdot3^k+228}{3^{k+6}}$, $\frac{608\cdot3^k+216}{3^{k+6}}$, \dots

\noindent $K[\omega^22+\omega5]=\frac{608}{729}$,\hspace{2.5cm} $K\cap(\frac{16384}{19683},\frac{608}{729})=[K[\omega^22+\omega5+1],K[\omega^22+\omega6])=[T_0[\omega^2 2+\omega5], T_0[\omega^2 2+\omega6])$\\
\noindent  $\frac{147712}{177147}$, $\frac{1329152}{1594323}$, $\frac{11960320}{14348907}$, $\frac{49216}{59049}$, $\frac{442880}{531441}$, $\frac{3985408}{4782969}$, $\frac{16400}{19683}$, $\frac{147584}{177147}$, $\frac{1328128}{1594323}$, $\frac{11952128}{14348907}$, $\frac{49184}{59049}$, $\frac{442624}{531441}$, $\frac{3983360}{4782969}$, $\frac{5464}{6561}$, $\frac{35848192}{43046721}$, $\frac{147520}{177147}$, $\frac{1327616}{1594323}$, $\frac{11948032}{14348907}$, $\frac{49168}{59049}$, $\frac{442496}{531441}$, $\frac{3982336}{4782969}$, $\frac{16388}{19683}$, $\frac{35840000}{43046721}$, $\frac{147488}{177147}$, $\frac{1327360}{1594323}$, $\frac{11945984}{14348907}$, $\frac{49160}{59049}$, $\frac{107511808}{129140163}$, $\frac{442432}{531441}$, $\frac{3981824}{4782969}$, $\frac{5462}{6561}$, $\frac{35835904}{43046721}$, $\frac{147472}{177147}$, $\frac{1327232}{1594323}$, $\frac{11944960}{14348907}$, $\frac{49156}{59049}$, $\frac{107503616}{129140163}$, $\frac{442400}{531441}$, $\frac{3981568}{4782969}$, $\frac{16385}{19683}$, $\frac{35833856}{43046721}$, $\frac{147464}{177147}$,\dots \\ $\frac{16384\cdot3^k+16384}{3^{k+9}}$, $\frac{16384\cdot3^k+15552}{3^{k+9}}$, $\frac{16384\cdot3^k+13824}{3^{k+9}}$, $\frac{16384\cdot3^k+13122}{3^{k+9}}$, $\frac{16384\cdot3^k+12288}{3^{k+9}}$, $\frac{16384\cdot3^k+11664}{3^{k+9}}$, $\frac{16384\cdot3^k+10368}{3^{k+9}}$,\\ $\frac{16384\cdot3^k+9216}{3^{k+9}}$, $\frac{16384\cdot3^k+8748}{3^{k+9}}$, $\frac{16384\cdot3^k+8192}{3^{k+9}}$, $\frac{16384\cdot3^k+7776}{3^{k+9}}$, $\frac{16384\cdot3^k+6912}{3^{k+9}}$, $\frac{16384\cdot3^k+6561}{3^{k+9}}$, $\frac{16384\cdot3^k+6144}{3^{k+9}}$, $\frac{16384\cdot3^k+5832}{3^{k+9}}$, \dots

\noindent $K[\omega^22+\omega6]=\frac{16384}{19683}$,\hspace{2.5cm} $K\cap(\frac{200}{243},\frac{16384}{19683})=[K[\omega^22+\omega6+1],K[\omega^22+\omega7])=[T_0[\omega^2 2+\omega6], T_0[\omega^2 2+\omega7])$\\
\noindent  $\frac{1820}{2187}$, $\frac{202}{243}$, $\mathbf{\frac{16352}{19683}}$, $\frac{5450}{6561}$, $\frac{605}{729}$, $\mathbf{\frac{146944}{177147}}$, $\frac{5440}{6561}$, $\mathbf{\frac{1813}{2187}}$, $\frac{604}{729}$, $\frac{16300}{19683}$, $\mathbf{\frac{5432}{6561}}$, $\frac{1810}{2187}$, $\frac{67}{81}$, $\mathbf{\frac{16280}{19683}}$, $\mathbf{\frac{48832}{59049}}$, $\frac{5425}{6561}$, $\frac{1808}{2187}$, $\frac{48800}{59049}$, $\frac{5420}{6561}$, $\mathbf{\frac{16256}{19683}}$, $\frac{602}{729}$, $\frac{16250}{19683}$, $\frac{1805}{2187}$, $\frac{16240}{19683}$, $\frac{1804}{2187}$, $\frac{48700}{59049}$, $\frac{5410}{6561}$, $\mathbf{\frac{146048}{177147}}$, $\frac{601}{729}$, $\frac{16225}{19683}$, $\frac{5408}{6561}$, $\frac{146000}{177147}$, $\frac{16220}{19683}$, $\frac{1802}{2187}$, $\frac{48650}{59049}$, $\frac{5405}{6561}$, $\frac{48640}{59049}$, $\frac{5404}{6561}$, $\frac{145900}{177147}$, $\frac{16210}{19683}$, $\frac{1801}{2187}$, $\frac{48625}{59049}$, $\frac{16208}{19683}$, $\frac{437600}{531441}$, $\frac{48620}{59049}$, $\frac{5402}{6561}$, $\frac{145850}{177147}$, $\frac{16205}{19683}$, $\frac{145840}{177147}$, $\frac{16204}{19683}$, $\frac{437500}{531441}$, $\frac{48610}{59049}$, $\frac{5401}{6561}$, $\frac{145825}{177147}$, $\frac{48608}{59049}$, \dots, $\frac{200\cdot3^k+200}{3^{k+5}}$, $\frac{200\cdot3^k+180}{3^{k+5}}$, $\frac{200\cdot3^k+162}{3^{k+5}}$, $\frac{200\cdot3^k+150}{3^{k+5}}$, $\frac{200\cdot3^k+135}{3^{k+5}}$, $\frac{200\cdot3^k+120}{3^{k+5}}$, $\frac{200\cdot3^k+108}{3^{k+5}}$, $\frac{200\cdot3^k+100}{3^{k+5}}$, $\frac{200\cdot3^k+90}{3^{k+5}}$, $\frac{200\cdot3^k+81}{3^{k+5}}$, $\frac{200\cdot3^k+75}{3^{k+5}}$, $\frac{200\cdot3^k+72}{3^{k+5}}$, \dots

\dots

\noindent $K[\omega^25+\omega3]=\frac{532}{729}$,\\ $K\cap(\frac{14336}{19683},\frac{532}{729})=[K[\omega^25+\omega3+1],K[\omega^25+\omega4]])=[T_0[\omega^25+\omega3], T_0[\omega^25+\omega4])$\\
\noindent  $\frac{1163264}{1594323}$, $\frac{129248}{177147}$, $\frac{1163008}{1594323}$, $\frac{43072}{59049}$, $\mathbf{\frac{14356}{19683}}$, $\frac{10465280}{14348907}$, $\mathbf{\frac{1595}{2187}}$, $\frac{387584}{531441}$, $\frac{43064}{59049}$, $\frac{387520}{531441}$, $\frac{4784}{6561}$, $\frac{3487232}{4782969}$, $\frac{129152}{177147}$, $\frac{14350}{19683}$, $\frac{1162240}{1594323}$, $\frac{129136}{177147}$, $\frac{1162112}{1594323}$, $\frac{43040}{59049}$, $\frac{10458112}{14348907}$, $\frac{387328}{531441}$, $\frac{43036}{59049}$, $\frac{3485696}{4782969}$, $\frac{387296}{531441}$, $\frac{14344}{19683}$, $\frac{3485440}{4782969}$, $\frac{129088}{177147}$, $\frac{4781}{6561}$, $\frac{1161728}{1594323}$, $\frac{129080}{177147}$, $\frac{1161664}{1594323}$, $\frac{43024}{59049}$, $\frac{10454528}{14348907}$, $\frac{387200}{531441}$, $\frac{43022}{59049}$, $\frac{3484672}{4782969}$, $\frac{387184}{531441}$, $\frac{4780}{6561}$, $\frac{3484544}{4782969}$, $\frac{129056}{177147}$, $\frac{1161472}{1594323}$, $\frac{129052}{177147}$, $\frac{10452992}{14348907}$, $\frac{1161440}{1594323}$, $\frac{43016}{59049}$, $\frac{10452736}{14348907}$, $\frac{387136}{531441}$, $\frac{43015}{59049}$, $\frac{3484160}{4782969}$, $\frac{387128}{531441}$, $\frac{14338}{19683}$, $\frac{3484096}{4782969}$, $\frac{129040}{177147}$, $\frac{1161344}{1594323}$, $\frac{129038}{177147}$, $\frac{10451968}{14348907}$, $\frac{1161328}{1594323}$, $\frac{43012}{59049}$, $\frac{10451840}{14348907}$, $\frac{387104}{531441}$, $\frac{3483904}{4782969}$, $\frac{387100}{531441}$, $\frac{59}{81}$, $\frac{3483872}{4782969}$, $\frac{129032}{177147}$, $\frac{1161280}{1594323}$, $\frac{129031}{177147}$, $\frac{10451456}{14348907}$, $\frac{1161272}{1594323}$, $\frac{43010}{59049}$, $\frac{10451392}{14348907}$, $\frac{387088}{531441}$, $\frac{3483776}{4782969}$, $\frac{387086}{531441}$, $\frac{3483760}{4782969}$, $\frac{129028}{177147}$, $\frac{1161248}{1594323}$, $\frac{10451200}{14348907}$, $\frac{1161244}{1594323}$, $\frac{43009}{59049}$, $\frac{10451168}{14348907}$, $\frac{387080}{531441}$, $\frac{3483712}{4782969}$, $\frac{387079}{531441}$, $\frac{3483704}{4782969}$, $\frac{129026}{177147}$, $\frac{1161232}{1594323}$, $\frac{10451072}{14348907}$, $\frac{1161230}{1594323}$, $\frac{10451056}{14348907}$, $\frac{387076}{531441}$, $\frac{3483680}{4782969}$, $\frac{3483676}{4782969}$, $\frac{129025}{177147}$, $\frac{1161224}{1594323}$, $\frac{10451008}{14348907}$, $\frac{1161223}{1594323}$, $\frac{10451000}{14348907}$, $\frac{387074}{531441}$, $\frac{3483664}{4782969}$, $\frac{3483662}{4782969}$, $\frac{1161220}{1594323}$, $\frac{10450976}{14348907}$, $\frac{10450972}{14348907}$, $\frac{387073}{531441}$, $\frac{3483656}{4782969}$, $\frac{3483655}{4782969}$,\dots, $\frac{14336\cdot3^k+13608}{3^{k+9}}$, $\frac{14336\cdot3^k+13122}{3^{k+9}}$, $\frac{14336\cdot3^k+11664}{3^{k+9}}$,  $\frac{14336\cdot3^k+10206}{3^{k+9}}$,  $\frac{14336\cdot3^k+8748}{3^{k+9}}$,  $\frac{14336\cdot3^k+7776}{3^{k+9}}$, $\frac{14336\cdot3^k+6804}{3^{k+9}}$, $\frac{14336\cdot3^k+6561}{3^{k+9}}$, $\frac{14336\cdot3^k+5832}{3^{k+9}}$,  $\frac{14336\cdot3^k+5103}{3^{k+9}}$, \dots

\dots

\normalsize
\section{Some elements of the compact \texorpdfstring{$K$: $T_1$}{K: T1}}

\tiny

\noindent $K[\omega]=K'[0]=\frac{4}{3}$,\hspace{3.5cm} $K'\cap(1,\frac{4}{3}]=[K'[0],K'[\omega])=[T_1[0], T_1[\omega])$\\
\noindent  $\frac{4}{3}$, $\mathbf{\frac{32}{27}}$, $\frac{10}{9}$, $\mathbf{\frac{256}{243}}$, $\frac{28}{27}$, $\frac{82}{81}$, $\frac{244}{243}$, $\frac{730}{729}$, $\frac{2188}{2187}$, $\frac{6562}{6561}$, $\frac{19684}{19683}$, $\frac{59050}{59049}$, $\frac{177148}{177147}$, $\frac{531442}{531441}$, $\frac{1594324}{1594323}$, $\frac{4782970}{4782969}$, $\frac{14348908}{14348907}$, \dots$\frac{3^k+1}{3^{k}}$, \dots

\noindent $K[\omega^2]=K'[\omega]=1$,\hspace{3.5cm} $K'\cap(\frac{8}{9},1)=(K'[\omega],K'[\omega2])=[T_1[\omega], T_1[\omega2])$\\
\noindent  $\frac{80}{81}$, $\frac{26}{27}$, $\frac{76}{81}$, $\mathbf{\frac{2048}{2187}}$, $\frac{25}{27}$, $\frac{224}{243}$, $\frac{74}{81}$, $\frac{220}{243}$, $\frac{73}{81}$, $\frac{656}{729}$, $\frac{218}{243}$, $\frac{652}{729}$, $\frac{217}{243}$, $\frac{1952}{2187}$, $\frac{650}{729}$, $\frac{1948}{2187}$, $\frac{649}{729}$, $\frac{5840}{6561}$, $\frac{1946}{2187}$, $\frac{5836}{6561}$, $\frac{1945}{2187}$, $\frac{17504}{19683}$, $\frac{5834}{6561}$, $\frac{17500}{19683}$, $\frac{5833}{6561}$, $\frac{52496}{59049}$, $\frac{17498}{19683}$, $\frac{52492}{59049}$, $\frac{17497}{19683}$, $\frac{157472}{177147}$, $\frac{52490}{59049}$, $\frac{157468}{177147}$, $\frac{52489}{59049}$,\dots$\frac{8\cdot3^k+8}{3^{k+2}}$, $\frac{8\cdot3^k+6}{3^{k+2}}$, $\frac{8\cdot3^k+4}{3^{k+2}}$, $\frac{8\cdot3^k+3}{3^{k+2}}$, \dots

\noindent $K[\omega^22]=K'[\omega2]=\frac{8}{9}$,\hspace{3.5cm} $K'\cap(\frac{64}{81},\frac{8}{9})=(K'[\omega2],K'[\omega3])=[T_1[\omega2], T_1[\omega3])$\\
\noindent  $\frac{640}{729}$, $\mathbf{\frac{70}{81}}$, $\frac{208}{243}$, $\frac{68}{81}$, $\frac{608}{729}$, $\mathbf{\frac{16384}{19683}}$, $\frac{200}{243}$, $\frac{1792}{2187}$, $\frac{22}{27}$, $\frac{592}{729}$, $\frac{196}{243}$, $\frac{1760}{2187}$, $\frac{65}{81}$, $\frac{584}{729}$, $\frac{5248}{6561}$, $\frac{194}{243}$, $\frac{1744}{2187}$, $\frac{580}{729}$, $\frac{5216}{6561}$, $\frac{193}{243}$, $\frac{1736}{2187}$, $\frac{15616}{19683}$, $\frac{578}{729}$, $\frac{5200}{6561}$, $\frac{1732}{2187}$, $\frac{15584}{19683}$, $\frac{577}{729}$, $\frac{5192}{6561}$, $\frac{46720}{59049}$, $\frac{1730}{2187}$, $\frac{15568}{19683}$, $\frac{5188}{6561}$, $\frac{46688}{59049}$, $\frac{1729}{2187}$, $\frac{15560}{19683}$,\dots$\frac{64\cdot3^k+64}{3^{k+4}}$, $\frac{64\cdot3^k+54}{3^{k+4}}$, $\frac{64\cdot3^k+48}{3^{k+4}}$, $\frac{64\cdot3^k+36}{3^{k+4}}$, $\frac{64\cdot3^k+32}{3^{k+4}}$, $\frac{64\cdot3^k+27}{3^{k+4}}$, $\frac{64\cdot3^k+24}{3^{k+4}}$, \dots

\noindent $K[\omega^23]=K'[\omega3]=\frac{64}{81}$,\hspace{3.5cm} $K'\cap(\frac{7}{9},\frac{64}{81})=(K'[\omega3],K'[\omega4])=[T_1[\omega3], T_1[\omega4])$\\
\noindent  $\frac{574}{729}$, $\frac{190}{243}$, $\frac{1708}{2187}$, $\mathbf{\frac{5120}{6561}}$, $\frac{568}{729}$, $\frac{5110}{6561}$, $\frac{1702}{2187}$, $\frac{15316}{19683}$, $\frac{5104}{6561}$, $\frac{45934}{59049}$, $\frac{15310}{19683}$, $\frac{137788}{177147}$, $\frac{45928}{59049}$, $\frac{413350}{531441}$, $\frac{137782}{177147}$, $\frac{1240036}{1594323}$, $\frac{413344}{531441}$, $\frac{3720094}{4782969}$, $\frac{1240030}{1594323}$, $\frac{11160268}{14348907}$, $\frac{3720088}{4782969}$, $\frac{33480790}{43046721}$, $\frac{11160262}{14348907}$, $\frac{100442356}{129140163}$, $\frac{33480784}{43046721}$, \dots$\frac{7\cdot3^k+7}{3^{k+2}}$, $\frac{7\cdot3^k+3}{3^{k+2}}$, \dots

\noindent $K[\omega^24]=K'[\omega4]=\frac{7}{9}$,\hspace{3.5cm} $K'\cap(\frac{20}{27},\frac{7}{9})=(K'[\omega4],K'[\omega5])=[T_1[\omega4], T_1[\omega5])$\\
\noindent  $\frac{560}{729}$, $\frac{62}{81}$, $\frac{185}{243}$, $\mathbf{\frac{1664}{2187}}$, $\frac{184}{243}$, $\frac{550}{729}$, $\frac{61}{81}$, $\frac{1640}{2187}$, $\frac{182}{243}$, $\frac{545}{729}$, $\frac{544}{729}$, $\frac{1630}{2187}$, $\frac{181}{243}$, $\frac{4880}{6561}$, $\frac{542}{729}$, $\frac{1625}{2187}$, $\frac{1624}{2187}$, $\frac{4870}{6561}$, $\frac{541}{729}$, $\frac{14600}{19683}$, $\frac{1622}{2187}$, $\frac{4865}{6561}$, $\frac{4864}{6561}$, $\frac{14590}{19683}$, $\frac{1621}{2187}$, $\frac{43760}{59049}$, $\frac{4862}{6561}$, $\frac{14585}{19683}$, $\frac{14584}{19683}$, $\frac{43750}{59049}$, $\frac{4861}{6561}$, \dots$\frac{20\cdot3^k+20}{3^{k+3}}$, $\frac{20\cdot3^k+18}{3^{k+3}}$, $\frac{20\cdot3^k+15}{3^{k+3}}$, $\frac{20\cdot3^k+12}{3^{k+3}}$, $\frac{20\cdot3^k+10}{3^{k+3}}$, $\frac{20\cdot3^k+9}{3^{k+3}}$, \dots

\noindent $K[\omega^25]=K'[\omega5]=\frac{20}{27}$,\hspace{3.5cm} $K'\cap(\frac{19}{27},\frac{20}{27})=(K'[\omega5],K'[\omega6])=[T_1[\omega5], T_1[\omega6])$\\
\noindent  $\mathbf{\frac{131072}{177147}}$, $\mathbf{\frac{1600}{2187}}$, $\frac{532}{729}$, $\mathbf{\frac{14336}{19683}}$, $\mathbf{\frac{176}{243}}$, $\mathbf{\frac{4736}{6561}}$, $\mathbf{\frac{175}{243}}$, $\mathbf{\frac{1568}{2187}}$, $\frac{58}{81}$, $\mathbf{\frac{14080}{19683}}$, $\mathbf{\frac{520}{729}}$, $\frac{1558}{2187}$, $\mathbf{\frac{4672}{6561}}$, $\mathbf{\frac{41984}{59049}}$, $\mathbf{\frac{518}{729}}$, $\mathbf{\frac{1552}{2187}}$, $\mathbf{\frac{13952}{19683}}$, $\frac{172}{243}$, $\mathbf{\frac{4640}{6561}}$, $\mathbf{\frac{41728}{59049}}$, $\frac{4636}{6561}$, $\mathbf{\frac{1544}{2187}}$, $\mathbf{\frac{13888}{19683}}$, $\mathbf{\frac{124928}{177147}}$, $\frac{514}{729}$, $\mathbf{\frac{4624}{6561}}$, $\frac{13870}{19683}$, $\mathbf{\frac{41600}{59049}}$, $\frac{1540}{2187}$, $\frac{41572}{59049}$, $\mathbf{\frac{13856}{19683}}$, $\frac{4618}{6561}$, $\frac{124678}{177147}$, $\mathbf{\frac{124672}{177147}}$, $\frac{13852}{19683}$, $\frac{373996}{531441}$, $\frac{41554}{59049}$, $\frac{1121950}{1594323}$, $\frac{124660}{177147}$, $\frac{3365812}{4782969}$, $\frac{373978}{531441}$, $\frac{10097398}{14348907}$, $\frac{1121932}{1594323}$, $\frac{30292156}{43046721}$, $\frac{3365794}{4782969}$, \dots$\frac{19\cdot3^k+19}{3^{k+3}}$, $\frac{19\cdot3^k+9}{3^{k+3}}$, \dots

\noindent $K[\omega^26]=K'[\omega6]=\frac{19}{27}$,\hspace{3.5cm} $K'\cap(\frac{512}{729},\frac{19}{27})=(K'[\omega6],K'[\omega7])=[T_1[\omega6], T_1[\omega7])$\\
\noindent  $\frac{4616}{6561}$, $\frac{41536}{59049}$, $\frac{373760}{531441}$, $\frac{1538}{2187}$, $\frac{13840}{19683}$, $\frac{124544}{177147}$, $\frac{4612}{6561}$, $\frac{41504}{59049}$, $\frac{373504}{531441}$, $\frac{1537}{2187}$, $\frac{13832}{19683}$, $\frac{124480}{177147}$, $\frac{1120256}{1594323}$, $\frac{4610}{6561}$, $\frac{41488}{59049}$, $\frac{373376}{531441}$, $\frac{13828}{19683}$, $\frac{124448}{177147}$, $\frac{1120000}{1594323}$, $\frac{4609}{6561}$, $\frac{41480}{59049}$, $\frac{373312}{531441}$, \dots$\frac{512\cdot3^k+512}{3^{k+6}}$, $\frac{512\cdot3^k+486}{3^{k+6}}$, $\frac{512\cdot3^k+432}{3^{k+6}}$, $\frac{512\cdot3^k+384}{3^{k+6}}$, $\frac{512\cdot3^k+324}{3^{k+6}}$, $\frac{512\cdot3^k+288}{3^{k+6}}$,\newline $\frac{512\cdot3^k+256}{3^{k+6}}$, $\frac{512\cdot3^k+243}{3^{k+6}}$, $\frac{512\cdot3^k+216}{3^{k+6}}$, $\frac{512\cdot3^k+192}{3^{k+6}}$, \dots

\noindent $K[\omega^27]=K'[\omega7]=\frac{512}{729}$,\hspace{3.5cm} $K'\cap(\frac{56}{81},\frac{512}{729})=(K'[\omega7],K'[\omega8])=[T_1[\omega7], T_1[\omega8])$\\
\noindent  $\frac{511}{729}$, $\frac{4592}{6561}$, $\frac{170}{243}$, $\frac{1526}{2187}$, $\frac{508}{729}$, $\frac{4564}{6561}$, $\frac{169}{243}$, $\frac{1520}{2187}$, $\frac{1519}{2187}$, $\frac{13664}{19683}$, $\frac{506}{729}$, $\mathbf{\frac{40960}{59049}}$, $\frac{4550}{6561}$, $\frac{1516}{2187}$, $\frac{13636}{19683}$, $\frac{505}{729}$, $\frac{4544}{6561}$, $\frac{4543}{6561}$, $\frac{40880}{59049}$, $\frac{1514}{2187}$, $\frac{13622}{19683}$, $\frac{4540}{6561}$, $\frac{40852}{59049}$, $\frac{1513}{2187}$, $\frac{13616}{19683}$, $\frac{13615}{19683}$, $\frac{122528}{177147}$, $\frac{4538}{6561}$, $\frac{40838}{59049}$, $\frac{13612}{19683}$, $\frac{122500}{177147}$, $\frac{4537}{6561}$, $\frac{40832}{59049}$, $\frac{40831}{59049}$, \dots$\frac{56\cdot3^k+56}{3^{k+4}}$, $\frac{56\cdot3^k+54}{3^{k+4}}$, $\frac{56\cdot3^k+42}{3^{k+4}}$, $\frac{56\cdot3^k+36}{3^{k+4}}$, $\frac{56\cdot3^k+28}{3^{k+4}}$, $\frac{56\cdot3^k+27}{3^{k+4}}$, $\frac{56\cdot3^k+24}{3^{k+4}}$, $\frac{56\cdot3^k+21}{3^{k+4}}$, \dots

\dots

\noindent $K[\omega^3]=K'[\omega^2]=\frac{2}{3}$,\hspace{3.5cm} $K'\cap(\frac{160}{243},\frac{2}{3})=(K'[\omega^2],K'[\omega^2+\omega])=[T_1[\omega^2], T_1[\omega^2+\omega])$\\
\noindent  $\frac{13120}{19683}$, $\frac{1456}{2187}$, $\frac{485}{729}$, $\frac{4360}{6561}$, $\frac{484}{729}$, $\frac{4352}{6561}$, $\frac{1450}{2187}$, $\frac{161}{243}$, $\frac{13040}{19683}$, $\frac{1448}{2187}$, $\frac{4340}{6561}$, $\frac{482}{729}$, $\frac{39040}{59049}$, $\frac{4336}{6561}$, $\frac{1445}{2187}$, $\frac{13000}{19683}$, $\frac{1444}{2187}$, $\frac{12992}{19683}$, $\frac{4330}{6561}$, $\frac{481}{729}$, $\frac{38960}{59049}$, $\frac{4328}{6561}$, $\frac{12980}{19683}$, $\frac{1442}{2187}$, \dots$\frac{160\cdot3^k+160}{3^{k+5}}$, $\frac{160\cdot3^k+144}{3^{k+5}}$, $\frac{160\cdot3^k+135}{3^{k+5}}$, $\frac{160\cdot3^k+120}{3^{k+5}}$, $\frac{160\cdot3^k+108}{3^{k+5}}$, $\frac{160\cdot3^k+96}{3^{k+5}}$, $\frac{160\cdot3^k+90}{3^{k+5}}$, $\frac{160\cdot3^k+81}{3^{k+5}}$, $\frac{160\cdot3^k+80}{3^{k+5}}$, $\frac{160\cdot3^k+72}{3^{k+5}}$, $\frac{160\cdot3^k+60}{3^{k+5}}$, $\frac{160\cdot3^k+54}{3^{k+5}}$, \dots

\noindent $K[\omega^3+\omega^2]=K'[\omega^2+\omega]=\frac{160}{243}$,\hspace{2cm} $K'\cap(\frac{52}{81},\frac{160}{243})=(K'[\omega^2+\omega],K'[\omega^2+\omega2])=[T_1[\omega^2+\omega], T_1[\omega^2+\omega2])$\\
\noindent  $\mathbf{\frac{1048576}{1594323}}$, $\frac{53}{81}$, $\frac{1430}{2187}$, $\mathbf{\frac{476}{729}}$, $\mathbf{\frac{475}{729}}$, $\mathbf{\frac{12800}{19683}}$, $\frac{158}{243}$, $\frac{4264}{6561}$, $\mathbf{\frac{4256}{6561}}$, $\frac{1417}{2187}$, $\frac{472}{729}$, $\mathbf{\frac{114688}{177147}}$, $\frac{157}{243}$, $\frac{4238}{6561}$, $\frac{470}{729}$, $\frac{12688}{19683}$, $\frac{4225}{6561}$, $\frac{1408}{2187}$, $\frac{469}{729}$, $\frac{12662}{19683}$, $\frac{1406}{2187}$, $\frac{37960}{59049}$, $\frac{12649}{19683}$, $\frac{4216}{6561}$, $\frac{1405}{2187}$, $\frac{37934}{59049}$, $\frac{4214}{6561}$, $\frac{113776}{177147}$, $\frac{37921}{59049}$, $\frac{12640}{19683}$, $\frac{4213}{6561}$, $\frac{113750}{177147}$, $\frac{12638}{19683}$, $\frac{341224}{531441}$, $\frac{113737}{177147}$, $\frac{37912}{59049}$, $\frac{12637}{19683}$, $\frac{341198}{531441}$, $\frac{37910}{59049}$, \dots$\frac{52\cdot3^k+52}{3^{k+4}}$, $\frac{52\cdot3^k+39}{3^{k+4}}$, $\frac{52\cdot3^k+36}{3^{k+4}}$, $\frac{52\cdot3^k+27}{3^{k+4}}$, $\frac{52\cdot3^k+26}{3^{k+4}}$, $\frac{52\cdot3^k+18}{3^{k+4}}$, \dots

\noindent $K[\omega^3+\omega^22]=K'[\omega^2+\omega2]=\frac{52}{81}$,\hspace{2cm} $K'\cap(\frac{17}{27},\frac{52}{81})=(K'[\omega^2+\omega2],K'[\omega^2+\omega3])=[T_1[\omega^2+\omega2], T_1[\omega^2+\omega3])$\\
\noindent  $\mathbf{\frac{37888}{59049}}$, $\mathbf{\frac{1400}{2187}}$, $\mathbf{\frac{155}{243}}$, $\frac{1394}{2187}$, $\mathbf{\frac{12544}{19683}}$, $\mathbf{\frac{4180}{6561}}$, $\mathbf{\frac{464}{729}}$, $\mathbf{\frac{112640}{177147}}$, $\mathbf{\frac{1387}{2187}}$, $\mathbf{\frac{4160}{6561}}$, $\frac{154}{243}$, $\mathbf{\frac{12464}{19683}}$, $\mathbf{\frac{37376}{59049}}$, $\frac{4148}{6561}$, $\mathbf{\frac{335872}{531441}}$, $\mathbf{\frac{4144}{6561}}$, $\mathbf{\frac{4142}{6561}}$, $\frac{460}{729}$, $\mathbf{\frac{12416}{19683}}$, $\frac{12410}{19683}$, $\frac{1378}{2187}$, $\mathbf{\frac{111616}{177147}}$, $\frac{37196}{59049}$, $\frac{4132}{6561}$, $\frac{111554}{177147}$, $\frac{12394}{19683}$, $\frac{334628}{531441}$, $\frac{37180}{59049}$, $\frac{1003850}{1594323}$, $\frac{111538}{177147}$, $\frac{3011516}{4782969}$, $\frac{334612}{531441}$, $\frac{9034514}{14348907}$, $\frac{1003834}{1594323}$, $\frac{27103508}{43046721}$, $\frac{3011500}{4782969}$, $\frac{81310490}{129140163}$, $\frac{9034498}{14348907}$, \dots$\frac{17\cdot3^k+17}{3^{k+3}}$, $\frac{17\cdot3^k+9}{3^{k+3}}$, \dots

\noindent $K[\omega^3+\omega^23]=K'[\omega^2+\omega3]=\frac{17}{27}$,\hspace{2cm} $K'\cap(\frac{152}{243},\frac{17}{27})=(K'[\omega^2+\omega3],K'[\omega^2+\omega4])=[T_1[\omega^2+\omega3], T_1[\omega^2+\omega4])$\\
\noindent  $\frac{12388}{19683}$, $\frac{1376}{2187}$, $\mathbf{\frac{1375}{2187}}$, $\mathbf{\frac{37120}{59049}}$, $\frac{4123}{6561}$, $\frac{458}{729}$, $\mathbf{\frac{333824}{531441}}$, $\frac{37088}{59049}$, $\mathbf{\frac{12352}{19683}}$, $\frac{12350}{19683}$, $\frac{1372}{2187}$, $\mathbf{\frac{111104}{177147}}$, $\frac{457}{729}$, $\mathbf{\frac{999424}{1594323}}$, $\frac{37012}{59049}$, $\frac{4112}{6561}$, $\frac{12331}{19683}$, $\mathbf{\frac{36992}{59049}}$, $\frac{1370}{2187}$, $\frac{110960}{177147}$, $\mathbf{\frac{332800}{531441}}$, $\frac{36974}{59049}$, $\frac{4108}{6561}$, $\frac{1369}{2187}$, $\frac{110884}{177147}$, $\frac{12320}{19683}$, $\frac{36955}{59049}$, $\frac{4106}{6561}$, $\frac{332576}{531441}$, $\mathbf{\frac{110848}{177147}}$, $\frac{110846}{177147}$, $\frac{12316}{19683}$, $\frac{4105}{6561}$, $\frac{332500}{531441}$, $\frac{36944}{59049}$, $\frac{110827}{177147}$, $\frac{12314}{19683}$, $\frac{997424}{1594323}$, $\frac{332462}{531441}$, $\frac{36940}{59049}$, $\mathbf{\frac{997376}{1594323}}$, $\frac{12313}{19683}$, $\frac{997348}{1594323}$, $\frac{110816}{177147}$, $\frac{332443}{531441}$, $\frac{36938}{59049}$, $\frac{2991968}{4782969}$, $\frac{997310}{1594323}$, $\frac{110812}{177147}$, $\frac{36937}{59049}$, $\frac{2991892}{4782969}$, $\frac{332432}{531441}$, $\frac{997291}{1594323}$, $\frac{110810}{177147}$, $\frac{8975600}{14348907}$, $\frac{2991854}{4782969}$, $\frac{332428}{531441}$, $\frac{110809}{177147}$, $\frac{8975524}{14348907}$, $\frac{997280}{1594323}$, $\frac{2991835}{4782969}$, $\frac{332426}{531441}$, \dots$\frac{152\cdot3^k+152}{3^{k+5}}$, $\frac{152\cdot3^k+114}{3^{k+5}}$, $\frac{152\cdot3^k+108}{3^{k+5}}$, $\frac{152\cdot3^k+81}{3^{k+5}}$, $\frac{152\cdot3^k+76}{3^{k+5}}$, $\frac{152\cdot3^k+72}{3^{k+5}}$, $\frac{152\cdot3^k+57}{3^{k+5}}$, $\frac{152\cdot3^k+54}{3^{k+5}}$, \dots

\noindent $K[\omega^3+\omega^24]=K'[\omega^2+\omega4]=\frac{152}{243}$,\hspace{2cm} $K'\cap(\frac{4096}{6561},\frac{152}{243})=(K'[\omega^2+\omega4],K'[\omega^2+\omega5])=[T_1[\omega^2+\omega4], T_1[\omega^2+\omega5])$\\
\noindent  $\frac{36928}{59049}$, $\frac{332288}{531441}$, $\frac{2990080}{4782969}$, $\frac{12304}{19683}$, $\frac{110720}{177147}$, $\frac{996352}{1594323}$, $\frac{4100}{6561}$, $\frac{36896}{59049}$, $\frac{332032}{531441}$, $\frac{2988032}{4782969}$, $\frac{12296}{19683}$, $\frac{110656}{177147}$, $\frac{995840}{1594323}$, $\frac{1366}{2187}$, $\frac{8962048}{14348907}$, $\frac{36880}{59049}$, $\frac{331904}{531441}$, $\frac{2987008}{4782969}$, $\frac{12292}{19683}$, $\frac{110624}{177147}$, $\frac{995584}{1594323}$, $\frac{4097}{6561}$, $\frac{8960000}{14348907}$, $\frac{36872}{59049}$, $\frac{331840}{531441}$, $\frac{2986496}{4782969}$, $\frac{12290}{19683}$, \dots$\frac{4096\cdot3^k+4096}{3^{k+8}}$, $\frac{4096\cdot3^k+3888}{3^{k+8}}$, $\frac{4096\cdot3^k+3456}{3^{k+8}}$, $\frac{4096\cdot3^k+3072}{3^{k+8}}$, $\frac{4096\cdot3^k+2916}{3^{k+8}}$, $\frac{4096\cdot3^k+2592}{3^{k+8}}$, $\frac{4096\cdot3^k+2304}{3^{k+8}}$, $\frac{4096\cdot3^k+2187}{3^{k+8}}$, $\frac{4096\cdot3^k+2048}{3^{k+8}}$, $\frac{4096\cdot3^k+1944}{3^{k+8}}$,\\ $\frac{4096\cdot3^k+1728}{3^{k+8}}$, $\frac{4096\cdot3^k+1536}{3^{k+8}}$, $\frac{4096\cdot3^k+1458}{3^{k+8}}$, \dots

\dots

\noindent $K[\omega^32]=K'[\omega^22]=\frac{16}{27}$,\hspace{2cm} $K'\cap(\frac{1280}{2187},\frac{16}{27})=(K'[\omega^22],K'[\omega^22+\omega])=[T_1[\omega^22], T_1[\omega^22+\omega])$\\
\noindent  $\frac{104960}{177147}$, $\mathbf{\frac{1295}{2187}}$, $\frac{11648}{19683}$, $\frac{3880}{6561}$, $\mathbf{\frac{1292}{2187}}$, $\frac{34880}{59049}$, $\frac{3872}{6561}$, $\frac{430}{729}$, $\frac{34816}{59049}$, $\frac{11600}{19683}$, $\frac{1288}{2187}$, $\frac{104320}{177147}$, $\frac{11584}{19683}$, $\mathbf{\frac{143}{243}}$, $\frac{3860}{6561}$, $\frac{34720}{59049}$, $\frac{3856}{6561}$, $\frac{312320}{531441}$, $\frac{1285}{2187}$, $\frac{34688}{59049}$, $\frac{11560}{19683}$, $\frac{428}{729}$, $\frac{104000}{177147}$, $\frac{11552}{19683}$, $\frac{3850}{6561}$, $\frac{103936}{177147}$, $\frac{34640}{59049}$, $\frac{3848}{6561}$, $\frac{311680}{531441}$, $\frac{34624}{59049}$, $\frac{11540}{19683}$, $\frac{1282}{2187}$, $\frac{103840}{177147}$, $\frac{11536}{19683}$, $\frac{934400}{1594323}$, $\frac{3845}{6561}$, $\frac{103808}{177147}$, $\frac{34600}{59049}$, $\frac{3844}{6561}$, $\frac{311360}{531441}$, $\frac{34592}{59049}$, $\frac{11530}{19683}$, $\frac{311296}{531441}$, $\frac{427}{729}$, $\frac{103760}{177147}$, $\frac{11528}{19683}$, \dots$\frac{1280\cdot3^k+1280}{3^{k+7}}$, $\frac{1280\cdot3^k+1215}{3^{k+7}}$, $\frac{1280\cdot3^k+1152}{3^{k+7}}$, $\frac{1280\cdot3^k+1080}{3^{k+7}}$, $\frac{1280\cdot3^k+972}{3^{k+7}}$, $\frac{1280\cdot3^k+960}{3^{k+7}}$, $\frac{1280\cdot3^k+864}{3^{k+7}}$, $\frac{1280\cdot3^k+810}{3^{k+7}}$, $\frac{1280\cdot3^k+768}{3^{k+7}}$, $\frac{1280\cdot3^k+729}{3^{k+7}}$, $\frac{1280\cdot3^k+720}{3^{k+7}}$, $\frac{1280\cdot3^k+648}{3^{k+7}}$, $\frac{1280\cdot3^k+640}{3^{k+7}}$, $\frac{1280\cdot3^k+576}{3^{k+7}}$, $\frac{1280\cdot3^k+540}{3^{k+7}}$, $\frac{1280\cdot3^k+486}{3^{k+7}}$, $\frac{1280\cdot3^k+480}{3^{k+7}}$, $\frac{1280\cdot3^k+432}{3^{k+7}}$, \dots

\noindent $K[\omega^32+\omega^2]=K'[\omega^22+\omega]=\frac{1280}{2187}$,\hspace{1cm} $K'\cap(\frac{140}{243},\frac{1280}{2187})=(K'[\omega^22+\omega],K'[\omega^22+\omega2])=[T_1[\omega^22+\omega], T_1[\omega^22+\omega2])$\\
\noindent  $\mathbf{\frac{8388608}{14348907}}$, $\frac{142}{243}$, $\frac{11480}{19683}$, $\frac{425}{729}$, $\frac{1274}{2187}$, $\frac{424}{729}$, $\frac{3815}{6561}$, $\mathbf{\frac{11440}{19683}}$, $\frac{1270}{2187}$, $\frac{3808}{6561}$, $\frac{47}{81}$, $\frac{11410}{19683}$, $\frac{1267}{2187}$, $\frac{3800}{6561}$, $\frac{422}{729}$, $\mathbf{\frac{3796}{6561}}$, $\frac{34160}{59049}$, $\frac{1265}{2187}$, $\frac{3794}{6561}$, $\mathbf{\frac{102400}{177147}}$, $\frac{1264}{2187}$, $\frac{11375}{19683}$, $\mathbf{\frac{34112}{59049}}$, $\frac{3790}{6561}$, $\frac{11368}{19683}$, $\frac{421}{729}$, $\frac{34090}{59049}$, $\frac{3787}{6561}$, $\frac{11360}{19683}$, $\frac{1262}{2187}$, $\frac{102200}{177147}$, $\frac{3785}{6561}$, $\frac{11354}{19683}$, $\frac{3784}{6561}$, $\frac{34055}{59049}$, $\frac{11350}{19683}$, $\frac{34048}{59049}$, $\frac{1261}{2187}$, $\frac{102130}{177147}$, $\frac{11347}{19683}$, $\frac{34040}{59049}$, $\frac{3782}{6561}$, \dots$\frac{140\cdot3^k+140}{3^{k+5}}$, $\frac{140\cdot3^k+135}{3^{k+5}}$, $\frac{140\cdot3^k+126}{3^{k+5}}$, $\frac{140\cdot3^k+108}{3^{k+5}}$, $\frac{140\cdot3^k+105}{3^{k+5}}$, $\frac{140\cdot3^k+90}{3^{k+5}}$, $\frac{140\cdot3^k+84}{3^{k+5}}$, $\frac{140\cdot3^k+81}{3^{k+5}}$, $\frac{140\cdot3^k+70}{3^{k+5}}$, $\frac{140\cdot3^k+63}{3^{k+5}}$, $\frac{140\cdot3^k+60}{3^{k+5}}$, $\frac{140\cdot3^k+54}{3^{k+5}}$, \dots

\normalsize
\section{Some elements of the compact \texorpdfstring{$K$: $T_2$}{K: T2}}

\tiny

\noindent $K[\omega^2]=K''[0]=1$,\hspace{3.5cm} $K''\cap(\frac{2}{3},1] =[K''[0],K''[\omega])=[T_2[0], T_2[\omega])$\\
\noindent  $1$, $\frac{8}{9}$, $\mathbf{\frac{64}{81}}$, $\frac{7}{9}$, $\frac{20}{27}$, $\frac{19}{27}$, $\mathbf{\frac{512}{729}}$, $\frac{56}{81}$, $\frac{55}{81}$, $\frac{164}{243}$, $\frac{163}{243}$, $\frac{488}{729}$, $\frac{487}{729}$, $\frac{1460}{2187}$, $\frac{1459}{2187}$, $\frac{4376}{6561}$, $\frac{4375}{6561}$, $\frac{13124}{19683}$, $\frac{13123}{19683}$,  \dots$\frac{2\cdot3^k+2}{3^{k+1}}$, $\frac{2\cdot3^k+1}{3^{k+1}}$, \dots

\noindent $K[\omega^3]=K''[\omega]=\frac{2}{3}$,\hspace{3.5cm} $K''\cap(\frac{16}{27},\frac{2}{3})=(K''[\omega],K''[\omega2])=[T_2[\omega], T_2[\omega2])$\\
\noindent  $\frac{160}{243}$, $\frac{52}{81}$, $\frac{17}{27}$, $\frac{152}{243}$, $\mathbf{\frac{4096}{6561}}$, $\frac{50}{81}$, $\frac{448}{729}$, $\frac{148}{243}$, $\frac{49}{81}$, $\frac{440}{729}$, $\frac{146}{243}$, $\frac{1312}{2187}$, $\frac{436}{729}$, $\frac{145}{243}$, $\frac{1304}{2187}$, $\frac{434}{729}$, $\frac{3904}{6561}$, $\frac{1300}{2187}$, $\frac{433}{729}$, $\frac{3896}{6561}$, $\frac{1298}{2187}$, $\frac{11680}{19683}$, $\frac{3892}{6561}$, $\frac{1297}{2187}$, $\frac{11672}{19683}$, $\frac{3890}{6561}$, $\frac{35008}{59049}$, $\frac{11668}{19683}$, $\frac{3889}{6561}$, $\frac{35000}{59049}$, $\frac{11666}{19683}$, \dots$\frac{16\cdot3^k+16}{3^{k+3}}$, $\frac{16\cdot3^k+12}{3^{k+3}}$, $\frac{16\cdot3^k+9}{3^{k+3}}$, $\frac{16\cdot3^k+8}{3^{k+3}}$, $\frac{16\cdot3^k+6}{3^{k+3}}$, \dots

\noindent $K[\omega^32]=K''[\omega2]=\frac{16}{27}$,\hspace{3.5cm} $K''\cap(\frac{5}{9},\frac{16}{27})=(K''[\omega2],K''[\omega3])=[T_2[\omega2], T_2[\omega3])$\\
\noindent  $\mathbf{\frac{1280}{2187}}$, $\frac{140}{243}$, $\mathbf{\frac{416}{729}}$, $\frac{46}{81}$, $\frac{410}{729}$, $\frac{136}{243}$, $\frac{1220}{2187}$, $\frac{406}{729}$, $\frac{3650}{6561}$, $\frac{1216}{2187}$, $\frac{10940}{19683}$, $\frac{3646}{6561}$, $\frac{32810}{59049}$, $\frac{10936}{19683}$, $\frac{98420}{177147}$, $\frac{32806}{59049}$, \dots$\frac{5\cdot3^k+5}{3^{k+2}}$, $\frac{5\cdot3^k+3}{3^{k+2}}$, \dots

\noindent $K[\omega^33]=K''[\omega3]=\frac{5}{9}$,\hspace{3.5cm} $K''\cap(\frac{128}{243},\frac{5}{9})=(K''[\omega3],K''[\omega4])=[T_2[\omega3], T_2[\omega4])$\\
\noindent  $\mathbf{\frac{32768}{59049}}$, $\frac{400}{729}$, $\mathbf{\frac{133}{243}}$, $\frac{3584}{6561}$, $\frac{44}{81}$, $\frac{1184}{2187}$, $\frac{392}{729}$, $\frac{3520}{6561}$, $\frac{130}{243}$, $\frac{1168}{2187}$, $\frac{10496}{19683}$, $\frac{388}{729}$, $\frac{3488}{6561}$, $\frac{43}{81}$, $\frac{1160}{2187}$, $\frac{10432}{19683}$, $\frac{386}{729}$, $\frac{3472}{6561}$, $\frac{31232}{59049}$, $\frac{1156}{2187}$, $\frac{10400}{19683}$, $\frac{385}{729}$, $\frac{3464}{6561}$, $\frac{31168}{59049}$, $\frac{1154}{2187}$, $\frac{10384}{19683}$, $\frac{93440}{177147}$, $\frac{3460}{6561}$, $\frac{31136}{59049}$, $\frac{1153}{2187}$, $\frac{10376}{19683}$, $\frac{93376}{177147}$, $\frac{3458}{6561}$, $\frac{31120}{59049}$,\dots$\frac{128\cdot3^k+128}{3^{k+5}}$, $\frac{128\cdot3^k+108}{3^{k+5}}$, $\frac{128\cdot3^k+96}{3^{k+5}}$, $\frac{128\cdot3^k+81}{3^{k+5}}$, $\frac{128\cdot3^k+72}{3^{k+5}}$, $\frac{128\cdot3^k+64}{3^{k+5}}$, $\frac{128\cdot3^k+54}{3^{k+5}}$, $\frac{128\cdot3^k+48}{3^{k+5}}$, \dots

\noindent $K[\omega^34]=K''[\omega4]=\frac{128}{243}$,\hspace{3.5cm} $K''\cap(\frac{14}{27},\frac{128}{243})=(K''[\omega4],K''[\omega5])=[T_2[\omega4], T_2[\omega5])$\\
\noindent  $\frac{1148}{2187}$, $\frac{127}{243}$, $\frac{1141}{2187}$, $\frac{380}{729}$, $\frac{3416}{6561}$, $\mathbf{\frac{10240}{19683}}$, $\frac{379}{729}$, $\frac{3409}{6561}$, $\frac{1136}{2187}$, $\frac{10220}{19683}$, $\frac{1135}{2187}$, $\frac{10213}{19683}$, $\frac{3404}{6561}$, $\frac{30632}{59049}$, $\frac{3403}{6561}$, $\frac{30625}{59049}$, $\frac{10208}{19683}$, $\frac{91868}{177147}$, $\frac{10207}{19683}$, $\frac{91861}{177147}$, $\frac{30620}{59049}$, $\frac{275576}{531441}$, $\frac{30619}{59049}$, $\frac{275569}{531441}$, $\frac{91856}{177147}$, $\frac{826700}{1594323}$, $\frac{91855}{177147}$, $\frac{826693}{1594323}$, $\frac{275564}{531441}$, \dots$\frac{14\cdot3^k+14}{3^{k+3}}$, $\frac{14\cdot3^k+9}{3^{k+3}}$, $\frac{14\cdot3^k+7}{3^{k+3}}$, $\frac{14\cdot3^k+6}{3^{k+3}}$, \dots

\noindent $K[\omega^35]=K''[\omega5]=\frac{14}{27}$,\hspace{3.5cm} $K''\cap(\frac{40}{81},\frac{14}{27})=(K''[\omega5],K''[\omega6])=[T_2[\omega5], T_2[\omega6])$\\
\noindent  $\frac{125}{243}$, $\frac{1120}{2187}$, $\frac{124}{243}$, $\frac{370}{729}$, $\mathbf{\frac{3328}{6561}}$, $\frac{41}{81}$, $\frac{368}{729}$, $\frac{1100}{2187}$, $\frac{122}{243}$, $\frac{365}{729}$, $\frac{3280}{6561}$, $\frac{364}{729}$, $\frac{1090}{2187}$, $\frac{121}{243}$, $\frac{1088}{2187}$, $\frac{3260}{6561}$, $\frac{362}{729}$, $\frac{1085}{2187}$, $\frac{9760}{19683}$, $\frac{1084}{2187}$, $\frac{3250}{6561}$, $\frac{361}{729}$, $\frac{3248}{6561}$, $\frac{9740}{19683}$, $\frac{1082}{2187}$, $\frac{3245}{6561}$, $\frac{29200}{59049}$, $\frac{3244}{6561}$, $\frac{9730}{19683}$, $\frac{1081}{2187}$, $\frac{9728}{19683}$, $\frac{29180}{59049}$, $\frac{3242}{6561}$, $\frac{9725}{19683}$,\dots$\frac{40\cdot3^k+40}{3^{k+4}}$, $\frac{40\cdot3^k+36}{3^{k+4}}$, $\frac{40\cdot3^k+30}{3^{k+4}}$, $\frac{40\cdot3^k+27}{3^{k+4}}$, $\frac{40\cdot3^k+24}{3^{k+4}}$, $\frac{40\cdot3^k+20}{3^{k+4}}$, $\frac{40\cdot3^k+18}{3^{k+4}}$, $\frac{40\cdot3^k+15}{3^{k+4}}$, \dots

\noindent $K[\omega^36]=K''[\omega6]=\frac{40}{81}$,\hspace{3.5cm} $K''\cap(\frac{13}{27},\frac{40}{81})=(K''[\omega6],K''[\omega7])=[T_2[\omega6], T_2[\omega7])$\\
\noindent  $\mathbf{\frac{262144}{531441}}$, $\mathbf{\frac{119}{243}}$, $\mathbf{\frac{3200}{6561}}$, $\frac{1066}{2187}$, $\mathbf{\frac{1064}{2187}}$, $\frac{118}{243}$, $\mathbf{\frac{28672}{59049}}$, $\frac{3172}{6561}$, $\frac{352}{729}$, $\frac{9490}{19683}$, $\frac{1054}{2187}$, $\frac{28444}{59049}$, $\frac{3160}{6561}$, $\frac{85306}{177147}$, $\frac{9478}{19683}$, $\frac{255892}{531441}$, $\frac{28432}{59049}$, $\frac{767650}{1594323}$, $\frac{85294}{177147}$, $\frac{2302924}{4782969}$, $\frac{255880}{531441}$, $\frac{6908746}{14348907}$, $\frac{767638}{1594323}$, $\frac{20726212}{43046721}$, $\frac{2302912}{4782969}$, $\frac{62178610}{129140163}$, $\frac{6908734}{14348907}$, \dots$\frac{13\cdot3^k+13}{3^{k+3}}$, $\frac{13\cdot3^k+9}{3^{k+3}}$, \dots

\noindent $K[\omega^37]=K''[\omega7]=\frac{13}{27}$,\hspace{3.5cm} $K''\cap(\frac{38}{81},\frac{13}{27})=(K''[\omega7],K''[\omega8])=[T_2[\omega7], T_2[\omega8])$\\
\noindent  $\mathbf{\frac{9472}{19683}}$, $\mathbf{\frac{350}{729}}$, $\mathbf{\frac{3136}{6561}}$, $\frac{1045}{2187}$, $\frac{116}{243}$, $\mathbf{\frac{28160}{59049}}$, $\mathbf{\frac{1040}{2187}}$, $\frac{3116}{6561}$, $\mathbf{\frac{9344}{19683}}$, $\mathbf{\frac{83968}{177147}}$, $\mathbf{\frac{1036}{2187}}$, $\frac{115}{243}$, $\mathbf{\frac{3104}{6561}}$, $\mathbf{\frac{27904}{59049}}$, $\frac{3097}{6561}$, $\frac{344}{729}$, $\mathbf{\frac{9280}{19683}}$, $\mathbf{\frac{83456}{177147}}$, $\frac{9272}{19683}$, $\mathbf{\frac{3088}{6561}}$, $\frac{343}{729}$, $\mathbf{\frac{27776}{59049}}$, $\mathbf{\frac{249856}{531441}}$, $\frac{9253}{19683}$, $\frac{1028}{2187}$, $\mathbf{\frac{9248}{19683}}$, $\frac{27740}{59049}$, $\mathbf{\frac{83200}{177147}}$, $\frac{1027}{2187}$, $\frac{27721}{59049}$, $\frac{3080}{6561}$, $\frac{83144}{177147}$, $\mathbf{\frac{27712}{59049}}$, $\frac{3079}{6561}$, $\frac{83125}{177147}$, $\frac{9236}{19683}$, $\frac{249356}{531441}$, $\frac{9235}{19683}$, $\mathbf{\frac{249344}{531441}}$, $\frac{249337}{531441}$, $\frac{27704}{59049}$, $\frac{747992}{1594323}$, $\frac{27703}{59049}$, $\frac{747973}{1594323}$, $\frac{83108}{177147}$, $\frac{2243900}{4782969}$, $\frac{83107}{177147}$, $\frac{2243881}{4782969}$, $\frac{249320}{531441}$, $\frac{6731624}{14348907}$, $\frac{249319}{531441}$, $\frac{6731605}{14348907}$, $\frac{747956}{1594323}$, \dots$\frac{38\cdot3^k+38}{3^{k+4}}$, $\frac{38\cdot3^k+27}{3^{k+4}}$, $\frac{38\cdot3^k+19}{3^{k+4}}$, $\frac{38\cdot3^k+18}{3^{k+4}}$, \dots

\noindent $K[\omega^38]=K''[\omega8]=\frac{38}{81}$,\hspace{3.5cm} $K''\cap(\frac{1024}{2187},\frac{38}{81})=(K''[\omega8],K''[\omega9])=[T_2[\omega8], T_2[\omega9])$\\
\noindent  $\frac{9232}{19683}$, $\frac{83072}{177147}$, $\frac{747520}{1594323}$, $\frac{3076}{6561}$, $\frac{27680}{59049}$, $\frac{249088}{531441}$, $\frac{1025}{2187}$, $\frac{9224}{19683}$, $\frac{83008}{177147}$, $\frac{747008}{1594323}$, $\frac{3074}{6561}$, $\frac{27664}{59049}$, $\frac{248960}{531441}$, $\frac{2240512}{4782969}$, $\frac{9220}{19683}$, $\frac{82976}{177147}$, $\frac{746752}{1594323}$, $\frac{3073}{6561}$, $\frac{27656}{59049}$, $\frac{248896}{531441}$, $\frac{2240000}{4782969}$, $\frac{9218}{19683}$, $\frac{82960}{177147}$, $\frac{746624}{1594323}$, \dots$\frac{1024\cdot3^k+1024}{3^{k+7}}$, $\frac{1024\cdot3^k+972}{3^{k+7}}$, $\frac{1024\cdot3^k+864}{3^{k+7}}$, $\frac{1024\cdot3^k+768}{3^{k+7}}$, $\frac{1024\cdot3^k+729}{3^{k+7}}$, $\frac{1024\cdot3^k+648}{3^{k+7}}$, $\frac{1024\cdot3^k+576}{3^{k+7}}$, $\frac{1024\cdot3^k+512}{3^{k+7}}$, $\frac{1024\cdot3^k+486}{3^{k+7}}$, $\frac{1024\cdot3^k+432}{3^{k+7}}$, $\frac{1024\cdot3^k+384}{3^{k+7}}$, \dots

\dots

\noindent $K[\omega^4]=K''[\omega^2]=\frac{4}{9}$,\hspace{3.5cm} $K''\cap(\frac{320}{729},\frac{4}{9})=(K''[\omega^2],K''[\omega^2+\omega])=[T_2[\omega^2], T_2[\omega^2+\omega])$\\
\noindent  $\frac{26240}{59049}$, $\frac{2912}{6561}$, $\frac{970}{2187}$, $\mathbf{\frac{323}{729}}$, $\frac{8720}{19683}$, $\frac{968}{2187}$, $\frac{8704}{19683}$, $\frac{2900}{6561}$, $\frac{322}{729}$, $\frac{26080}{59049}$, $\frac{2896}{6561}$, $\frac{965}{2187}$, $\frac{8680}{19683}$, $\frac{964}{2187}$, $\frac{78080}{177147}$, $\frac{8672}{19683}$, $\frac{2890}{6561}$, $\frac{107}{243}$, $\frac{26000}{59049}$, $\frac{2888}{6561}$, $\frac{25984}{59049}$, $\frac{8660}{19683}$, $\frac{962}{2187}$, $\frac{77920}{177147}$, $\frac{8656}{19683}$, $\frac{2885}{6561}$, $\frac{25960}{59049}$, $\frac{2884}{6561}$, $\frac{233600}{531441}$, $\frac{25952}{59049}$, $\frac{8650}{19683}$, $\frac{961}{2187}$, $\frac{77840}{177147}$,$\frac{8648}{19683}$, $\frac{77824}{177147}$,\dots$\frac{320\cdot3^k+320}{3^{k+6}}$, $\frac{320\cdot3^k+288}{3^{k+6}}$, $\frac{320\cdot3^k+270}{3^{k+6}}$, $\frac{320\cdot3^k+243}{3^{k+6}}$, $\frac{320\cdot3^k+240}{3^{k+6}}$, $\frac{320\cdot3^k+216}{3^{k+6}}$, $\frac{320\cdot3^k+192}{3^{k+6}}$, $\frac{320\cdot3^k+180}{3^{k+6}}$, $\frac{320\cdot3^k+162}{3^{k+6}}$, $\frac{320\cdot3^k+160}{3^{k+6}}$, $\frac{320\cdot3^k+144}{3^{k+6}}$, $\frac{320\cdot3^k+135}{3^{k+6}}$, $\frac{320\cdot3^k+120}{3^{k+6}}$, $\frac{320\cdot3^k+108}{3^{k+6}}$, \dots

\noindent $K[\omega^4+\omega^3]=K''[\omega^2+\omega]=\frac{320}{729}$,\hspace{2.5cm} $K''\cap(\frac{35}{81},\frac{320}{729})=(K''[\omega^2+\omega],K''[\omega^2+\omega2])=[T_2[\omega^2+\omega], T_2[\omega^2+\omega2])$\\
\noindent  $\mathbf{\frac{2097152}{4782969}}$, $\frac{2870}{6561}$, $\frac{106}{243}$, $\mathbf{\frac{2860}{6561}}$, $\frac{952}{2187}$, $\frac{950}{2187}$, $\mathbf{\frac{949}{2187}}$, $\frac{8540}{19683}$, $\mathbf{\frac{25600}{59049}}$, $\frac{316}{729}$, $\mathbf{\frac{8528}{19683}}$, $\frac{2842}{6561}$, $\frac{2840}{6561}$, $\frac{25550}{59049}$, $\frac{946}{2187}$, $\frac{8512}{19683}$, $\frac{8510}{19683}$, $\frac{76580}{177147}$, $\frac{2836}{6561}$, $\frac{25522}{59049}$, $\frac{25520}{59049}$, $\frac{229670}{531441}$, $\frac{8506}{19683}$, $\frac{76552}{177147}$, $\frac{76550}{177147}$,  \dots$\frac{35\cdot3^k+35}{3^{k+4}}$, $\frac{35\cdot3^k+27}{3^{k+4}}$, $\frac{35\cdot3^k+21}{3^{k+4}}$, $\frac{35\cdot3^k+15}{3^{k+4}}$, \dots

\noindent $K[\omega^4+\omega^32]=K''[\omega^2+\omega2]=\frac{35}{81}$,\hspace{2cm} $K''\cap(\frac{104}{243},\frac{35}{81})=(K''[\omega^2+\omega2],K''[\omega^2+\omega3])=[T_2[\omega^2+\omega2], T_2[\omega^2+\omega3])$\\
\noindent  $\frac{2834}{6561}$, $\frac{944}{2187}$, $\mathbf{\frac{229376}{531441}}$, $\frac{314}{729}$, $\frac{8476}{19683}$, $\frac{2821}{6561}$, $\frac{940}{2187}$, $\frac{25376}{59049}$, $\frac{313}{729}$, $\frac{8450}{19683}$, $\frac{2816}{6561}$, $\frac{938}{2187}$, $\frac{25324}{59049}$, $\frac{8437}{19683}$, $\frac{2812}{6561}$, $\frac{75920}{177147}$, $\frac{937}{2187}$, $\frac{25298}{59049}$, $\frac{8432}{19683}$, $\frac{2810}{6561}$, $\frac{75868}{177147}$, $\frac{25285}{59049}$, $\frac{8428}{19683}$, \dots$\frac{104\cdot3^k+104}{3^{k+5}}$, $\frac{104\cdot3^k+81}{3^{k+5}}$, $\frac{104\cdot3^k+78}{3^{k+5}}$, $\frac{104\cdot3^k+72}{3^{k+5}}$, $\frac{104\cdot3^k+54}{3^{k+5}}$, $\frac{104\cdot3^k+52}{3^{k+5}}$, $\frac{104\cdot3^k+39}{3^{k+5}}$, $\frac{104\cdot3^k+36}{3^{k+5}}$, \dots

\noindent $K[\omega^4+\omega^33]=K''[\omega^2+\omega3]=\frac{104}{243}$,\hspace{2cm} $K''\cap(\frac{34}{81},\frac{104}{243})=(K''[\omega^2+\omega3],K''[\omega^2+\omega4])=[T_2[\omega^2+\omega3], T_2[\omega^2+\omega4])$\\
\noindent  $\mathbf{\frac{75776}{177147}}$, $\frac{935}{2187}$, $\mathbf{\frac{2800}{6561}}$, $\mathbf{\frac{931}{2187}}$, $\mathbf{\frac{310}{729}}$, $\frac{2788}{6561}$, $\mathbf{\frac{25088}{59049}}$, $\mathbf{\frac{8360}{19683}}$, $\mathbf{\frac{928}{2187}}$, $\mathbf{\frac{225280}{531441}}$, $\frac{103}{243}$, $\mathbf{\frac{925}{2187}}$, $\mathbf{\frac{2774}{6561}}$, $\mathbf{\frac{8320}{19683}}$, $\frac{308}{729}$, $\frac{2771}{6561}$, $\mathbf{\frac{24928}{59049}}$, $\mathbf{\frac{74752}{177147}}$, $\frac{8296}{19683}$, $\mathbf{\frac{671744}{1594323}}$, $\frac{307}{729}$, $\mathbf{\frac{8288}{19683}}$, $\mathbf{\frac{8284}{19683}}$, $\frac{920}{2187}$, $\frac{8279}{19683}$, $\mathbf{\frac{24832}{59049}}$, $\frac{24820}{59049}$, $\frac{919}{2187}$, $\frac{2756}{6561}$, $\mathbf{\frac{223232}{531441}}$, $\frac{24803}{59049}$, $\frac{74392}{177147}$, $\frac{2755}{6561}$, $\frac{8264}{19683}$, $\frac{74375}{177147}$, $\frac{223108}{531441}$, $\frac{8263}{19683}$, $\frac{24788}{59049}$, $\frac{223091}{531441}$, \dots$\frac{34\cdot3^k+34}{3^{k+4}}$, $\frac{34\cdot3^k+27}{3^{k+4}}$, $\frac{34\cdot3^k+18}{3^{k+4}}$, $\frac{34\cdot3^k+17}{3^{k+4}}$, \dots

\noindent $K[\omega^4+\omega^34]=K''[\omega^2+\omega4]=\frac{34}{81}$,\hspace{2cm} $K''\cap(\frac{304}{729},\frac{34}{81})=(K''[\omega^2+\omega4],K''[\omega^2+\omega5])=[T_2[\omega^2+\omega4], T_2[\omega^2+\omega5])$\\
\noindent  $\frac{24776}{59049}$, $\frac{2752}{6561}$, $\mathbf{\frac{2750}{6561}}$, $\mathbf{\frac{74240}{177147}}$, $\frac{8246}{19683}$, $\frac{916}{2187}$, $\mathbf{\frac{667648}{1594323}}$, $\frac{74176}{177147}$, $\frac{305}{729}$, $\mathbf{\frac{24704}{59049}}$, $\frac{24700}{59049}$, $\frac{2744}{6561}$, $\mathbf{\frac{222208}{531441}}$, $\frac{8227}{19683}$, $\frac{914}{2187}$, $\mathbf{\frac{1998848}{4782969}}$, $\frac{74024}{177147}$, $\frac{8224}{19683}$, $\frac{24662}{59049}$, $\mathbf{\frac{73984}{177147}}$, $\frac{2740}{6561}$, $\frac{221920}{531441}$, $\mathbf{\frac{665600}{1594323}}$, $\frac{913}{2187}$, $\frac{73948}{177147}$, $\frac{8216}{19683}$, $\frac{24643}{59049}$, $\frac{2738}{6561}$, $\frac{221768}{531441}$, $\frac{24640}{59049}$, $\frac{73910}{177147}$, $\frac{8212}{19683}$, $\frac{665152}{1594323}$, $\frac{2737}{6561}$, $\mathbf{\frac{221696}{531441}}$, $\frac{221692}{531441}$, $\frac{24632}{59049}$, $\frac{73891}{177147}$, $\frac{8210}{19683}$, $\frac{665000}{1594323}$, $\frac{73888}{177147}$, $\frac{221654}{531441}$, $\frac{24628}{59049}$, $\frac{1994848}{4782969}$, $\frac{8209}{19683}$, $\frac{664924}{1594323}$, $\frac{73880}{177147}$, $\mathbf{\frac{1994752}{4782969}}$, $\frac{221635}{531441}$, $\frac{24626}{59049}$, $\frac{1994696}{4782969}$, $\frac{221632}{531441}$, $\frac{664886}{1594323}$, $\frac{73876}{177147}$, $\frac{5983936}{14348907}$, $\frac{24625}{59049}$, $\frac{1994620}{4782969}$, $\frac{221624}{531441}$, $\frac{664867}{1594323}$, $\frac{73874}{177147}$, $\frac{5983784}{14348907}$, $\frac{664864}{1594323}$, $\frac{1994582}{4782969}$, $\frac{221620}{531441}$, \dots$\frac{304\cdot3^k+304}{3^{k+6}}$, $\frac{304\cdot3^k+243}{3^{k+6}}$, $\frac{304\cdot3^k+228}{3^{k+6}}$, $\frac{304\cdot3^k+216}{3^{k+6}}$, $\frac{304\cdot3^k+171}{3^{k+6}}$, $\frac{304\cdot3^k+162}{3^{k+6}}$, $\frac{304\cdot3^k+152}{3^{k+6}}$, $\frac{304\cdot3^k+144}{3^{k+6}}$, $\frac{304\cdot3^k+114}{3^{k+6}}$, $\frac{304\cdot3^k+108}{3^{k+6}}$, \dots

\dots

\noindent $K[\omega^42]=K''[\omega^22]=\frac{32}{81}$,\hspace{2cm} $K''\cap(\frac{95}{243},\frac{32}{81})=(K''[\omega^22],K''[\omega^22+\omega])=[T_2[\omega^22], T_2[\omega^22+\omega])$\\
\noindent  $\mathbf{\frac{209920}{531441}}$, $\mathbf{\frac{2590}{6561}}$, $\mathbf{\frac{23296}{59049}}$, $\mathbf{\frac{7760}{19683}}$, $\frac{2584}{6561}$, $\mathbf{\frac{69760}{177147}}$, $\mathbf{\frac{287}{729}}$, $\mathbf{\frac{7744}{19683}}$, $\frac{860}{2187}$, $\mathbf{\frac{69632}{177147}}$, $\mathbf{\frac{23200}{59049}}$, $\mathbf{\frac{2576}{6561}}$, $\mathbf{\frac{208640}{531441}}$, $\frac{23180}{59049}$, $\mathbf{\frac{23168}{59049}}$, $\frac{286}{729}$, $\mathbf{\frac{7720}{19683}}$, $\mathbf{\frac{69440}{177147}}$, $\frac{7714}{19683}$, $\mathbf{\frac{7712}{19683}}$, $\mathbf{\frac{624640}{1594323}}$, $\frac{2570}{6561}$, $\mathbf{\frac{69376}{177147}}$, $\mathbf{\frac{23120}{59049}}$, $\frac{69350}{177147}$, $\frac{856}{2187}$, $\mathbf{\frac{208000}{531441}}$, $\frac{23104}{59049}$, $\frac{7700}{19683}$, $\mathbf{\frac{207872}{531441}}$, $\frac{207860}{531441}$, $\frac{2566}{6561}$, $\mathbf{\frac{69280}{177147}}$, $\frac{69274}{177147}$, $\frac{23090}{59049}$, $\frac{623390}{1594323}$, $\frac{7696}{19683}$, $\mathbf{\frac{623360}{1594323}}$, $\frac{207784}{531441}$, $\frac{69260}{177147}$, $\frac{1869980}{4782969}$, $\frac{23086}{59049}$, $\frac{623314}{1594323}$, $\frac{207770}{531441}$, $\frac{5609750}{14348907}$, $\frac{69256}{177147}$, $\frac{1869904}{4782969}$, $\frac{623300}{1594323}$, $\frac{16829060}{43046721}$, $\frac{207766}{531441}$, $\frac{5609674}{14348907}$, $\frac{1869890}{4782969}$, \dots$\frac{95\cdot3^k+95}{3^{k+5}}$, $\frac{95\cdot3^k+81}{3^{k+5}}$, $\frac{95\cdot3^k+57}{3^{k+5}}$, $\frac{95\cdot3^k+45}{3^{k+5}}$, \dots

\noindent $K[\omega^42+\omega^3]=K''[\omega^22+\omega]=\frac{95}{243}$,\hspace{1.5cm} $K''\cap(\frac{2560}{6561},\frac{95}{243})=(K''[\omega^22+\omega],K''[\omega^22+\omega2])=[T_2[\omega^22+\omega], T_2[\omega^22+\omega2])$\\
\noindent  $\frac{69248}{177147}$, $\frac{23080}{59049}$, $\frac{2564}{6561}$, $\frac{207680}{531441}$, $\frac{23072}{59049}$, $\frac{1868800}{4782969}$, $\frac{7690}{19683}$, $\frac{207616}{531441}$, $\frac{69200}{177147}$, $\frac{7688}{19683}$, $\frac{622720}{1594323}$, $\frac{69184}{177147}$, $\frac{23060}{59049}$, $\frac{622592}{1594323}$, $\frac{854}{2187}$, $\frac{207520}{531441}$, $\frac{23056}{59049}$, $\frac{1867520}{4782969}$, $\frac{7685}{19683}$, $\frac{207488}{531441}$, $\frac{69160}{177147}$, $\frac{7684}{19683}$, $\frac{622400}{1594323}$, $\frac{69152}{177147}$, $\frac{5601280}{14348907}$, $\frac{23050}{59049}$, $\frac{622336}{1594323}$, $\frac{2561}{6561}$, $\frac{207440}{531441}$, $\frac{23048}{59049}$, $\frac{1866880}{4782969}$, $\frac{207424}{531441}$, $\frac{69140}{177147}$, $\frac{1866752}{4782969}$, $\frac{7682}{19683}$, $\frac{622240}{1594323}$, $\frac{69136}{177147}$, $\frac{5600000}{14348907}$, $\frac{23045}{59049}$, $\frac{622208}{1594323}$, $\frac{207400}{531441}$, $\frac{23044}{59049}$, $\frac{1866560}{4782969}$, $\frac{207392}{531441}$,\dots \\ $\frac{2560\cdot3^k+2560}{3^{k+8}}$,\hfill $\frac{2560\cdot3^k+2430}{3^{k+8}}$,\hfill $\frac{2560\cdot3^k+2304}{3^{k+8}}$,\hfill $\frac{2560\cdot3^k+2187}{3^{k+8}}$,\hfill $\frac{2560\cdot3^k+2160}{3^{k+8}}$,\hfill $\frac{2560\cdot3^k+1944}{3^{k+8}}$,\hfill $\frac{2560\cdot3^k+1920}{3^{k+8}}$,\hfill $\frac{2560\cdot3^k+1728}{3^{k+8}}$,\\ $\frac{2560\cdot3^k+1620}{3^{k+8}}$,\hfill $\frac{2560\cdot3^k+1536}{3^{k+8}}$,\hfill $\frac{2560\cdot3^k+1458}{3^{k+8}}$,\hfill $\frac{2560\cdot3^k+1440}{3^{k+8}}$,\hfill $\frac{2560\cdot3^k+1296}{3^{k+8}}$,\hfill $\frac{2560\cdot3^k+1280}{3^{k+8}}$,\hfill $\frac{2560\cdot3^k+1215}{3^{k+8}}$,\hfill $\frac{2560\cdot3^k+1152}{3^{k+8}}$,\\ $\frac{2560\cdot3^k+1080}{3^{k+8}}$,$\frac{2560\cdot3^k+972}{3^{k+8}}$, $\frac{2560\cdot3^k+960}{3^{k+8}}$, $\frac{2560\cdot3^k+864}{3^{k+8}}$, \dots

} 

\end{document}